\documentclass[reqno,makeidx]{amsart}
        
\input supp-pdf.tex
\usepackage{graphicx}
\usepackage{latexsym}
\usepackage{amstext}
\usepackage {amsmath}
\usepackage {amsfonts}
\usepackage {amssymb}
\usepackage {amsthm}
\usepackage {bbm}
\usepackage{enumerate}
\usepackage{array}
\usepackage{comment}
\usepackage{xspace}
\usepackage[bookmarks=true]{hyperref}
\usepackage{stmaryrd}

\DeclareMathAlphabet{\mathpzc}{OT1}{pzc}{m}{it}

\newcommand{\res}      {\mathop{\hbox{\vrule height 7pt width .5pt depth
                        0pt\vrule height .5pt width 6pt depth
0pt}}\nolimits}

\newcommand{\C}{\mathcal{C}}
\def\calBepsu{B^\eps_u} 
 \def\Hepsu{H^\eps_u} 

\newcommand{\constensurecompact}{c}
\newcommand{\constant}{C}

\newcommand{\elbend}{\varkappa_b}
\newcommand{\elGauss}{\varkappa_G}

\newcommand{\eps}{\varepsilon}
\def\fepsu{{f^\eps_u}} 
\def\fepsueps{f^\eps_{u_\eps}} 
\newcommand{\gae}{\gamma_\eps}
\newcommand{\grad}{\nabla}
\newcommand{\K}{\ensuremath{\mathpzc K}}

\newcommand{\Ha}{\ensuremath{\mathcal{H}}}
\newcommand{\Helf}{\mathpzc W_{{\rm Hel}}}
\def\Hepsueps{H^\eps_{u_\eps}}

\newcommand{\Hv}{\mathbf{H}}

\newcommand{\LL}{\ensuremath{\mathcal{L}}}

\def\muepsu{\mu^\eps_u}
\def\muepsueps{\mu^\eps_{u_\eps}}
\newcommand{\N}{\ensuremath{\mathbb{N}}}
\def\niu{\nu_u}
\def\niueps{\nu_{u_\eps}}
\def\niuepsk{\nu_{u_{\eps_k}}}

\def\proju{P^u}
\def\projueps{P^{u_\eps}}
\newcommand{\q}{\mathbf{q}}
\def\R{\mathbb R}
\def\rectifiableset{\mathcal M}
\newcommand{\Rn}{\ensuremath{{\mathbb{R}^n}}}

\newcommand{\sffeps}{\mathbf{B}_{u_\eps}}
\newcommand{\sffu}{\mathbf{B}_{u}}
\newcommand{\sff}{\mathbf B}

\newcommand{\surftens}{c_0}

\newcommand{\tgae}{\widetilde{\gamma}_\eps}
\def\tildemuepsu{\widetilde\mu^\eps_u}
\def\tildemuepsueps{\widetilde\mu^\eps_{u_\eps}} 
\def\tildemuepsuepsk{\widetilde\mu^{\eps_k}_{u_{\eps_k}}} 
\newcommand{\var}{\mathbf v}
\def\Vepsu{V^\eps_u} 
\def\Vepsueps{V^\eps_{u_\eps}} 

\def\Vzeroepsu{V^{0,\eps}_u}
\def\Vzeroepsueps{V^{0,\eps}_{u_\eps}} 
\def\Vzeroepsuepsk{V^{0,{\eps_k}}_{u_{\eps_k}}} 
\def\xiepsu{\xi^\eps_u}
\def\xiepsueps{\xi^\eps_{u_\eps}} 
\def\xiepsuepsk{\xi^{\eps_k}_{u_{\eps_k}}} 
\def\BBu{{\bf B}_u}

\theoremstyle{plain}
\numberwithin{equation}{section}
\newtheorem{lemma}{Lemma}[section]
\newtheorem{theorem}[lemma]{Theorem}
\newtheorem{proposition}[lemma]{Proposition}
\newtheorem{definition}[lemma]{Definition}

\newtheorem{corollary}[lemma]{Corollary}
\theoremstyle{definition}
\newtheorem{remark}[lemma]{Remark}
\begin{document}
\title{Approximation of the Helfrich's functional via Diffuse Interfaces}
\author{Giovanni Bellettini}
\address{Giovanni Bellettini, Dipartimento di Matematica, 
Universit\`a di Roma Tor Vergata, via della Ricerca Scientifica,
00133 Roma, Italy, and Laboratori Nazionali di Frascati,
Istituto Nazionale di Fisica Nucleare (INFN)
via E. Fermi, 40, Frascati (Roma), I-00044, Italy}
\email{Giovanni.Bellettini@lnf.infn.it}

\author{Luca Mugnai}
\address{Luca Mugnai, Max Planck Institute for Mathematics in the
  Sciences, Inselstr. 22, D-04103 Leipzig, Germany}
\email{mugnai@mis.mpg.de, }



\date{\today}
\begin{abstract}
We give a rigorous proof of 
the approximability 
of the so-called Helfrich's 
functional via diffuse interfaces, 
under a constraint on the ratio between the bending rigidity 
and the Gauss-rigidity.
\end{abstract}

\maketitle
\section{Introduction}\label{sec:int}
Let $\Omega\subset\R^3$ be an open connected set with smooth 
boundary.  
Define
\begin{equation}\label{eq:Helfrich}
\Helf(E):=\int_{\Omega \cap \partial 
E}\left[\frac{\elbend}{2}(H_{\partial 
E}-H_0)^2+\elGauss K_{\partial 
E}\right]\,d\Ha^{2},
\end{equation}
where $E\subset \Omega$ is open, bounded and with boundary 
$\partial E$ of class $C^2$ in $\Omega$; 
 $H_{\partial E},\, K_{\partial E}$ are respectively the mean
curvature  and the Gaussian-curvature of $\partial E$ (i.e. respectively 
the sum and the product of the two principal curvatures of $\partial E$); 
$\Ha^2$ is the $2$-dimensional Hausdorff-measure; 
$\elbend,\,H_0,\,\elGauss$ are given constants.
For our purposes it is convenient to write $\Helf$ as
\begin{equation*}
\Helf(E)=\frac{\elbend}{2} \mathpzc H(E)+\elGauss \mathpzc 
K(E),
\end{equation*}
where
\begin{align}
&\mathpzc{H}(E):=\int_{\Omega\cap\partial E}\left(H_{\partial E}-H_0\right)^2\,d\Ha^2,
\label{comeunvecchiorimorsoounvizioassurdo}
\\
&\mathpzc{K}(E):=\int_{\Omega\cap\partial E} K_{\partial E}~d\Ha^2.
\label{arturoantunesdecoimbradettozicooilfolletto}
\end{align}

The functional $\Helf$ was proposed by Helfrich as a surface 
 energy for closed biological membranes represented by a smooth 
boundaryless surface (see also \cite{Can, Eva} and \cite[Chapter 
7]{Boal}).  Minimizers and critical points of $\Helf$ in the 
 class of subsets $E\subset\Omega$ satisfying a constraint on the  area 
 $\Ha^2(\Omega\cap\partial E)$ and on the enclosed volume 
 $\LL^3(E\cap\Omega)$, are expected to describe approximately  the shape 
 of biological membranes such as monolayers or lipid bilayers (see again 
 \cite{Boal} for an introduction to the subject).
Note that  the term $\mathpzc K(E)$ can be neglected when minimizing $\Helf(E)$ under a topological constraint 
on $E$, 
since by the Gauss-Bonnet theorem it 
reduces to a constant depending on the fixed topology.  On 
the other hand $\mathpzc K$ plays an essential role 
in several recent related models 
(see e.g. \cite{BaumgartNature, AmAll, ArrDeSim}).  

The constant $\elbend>0$ is called the bending rigidity.
 The constant $H_0$ is called the spontaneous 
 curvature. It  is expected to be non zero  
 when dealing with biological membranes such as bilayers  with chemically 
 different interior and exterior layers, or when different enviroments 
 inside and outside the membrane are source of asymmetry.
Observe  that,  when $H_0\neq 0$, the 
functional $\mathpzc H$ depends on the orientation of $\partial E$ (and not only on the geometry of $\partial E$ as 
in the case $H_0=0$). 
 The constant $\elGauss$ is called the Gauss-rigidity.
 Although few experimental measurements for $\elGauss$ are 
presently available,  it is expected to be negative (see \cite{TKS-Gauss}, 
 \cite{SiegKoz-Gauss}, \cite[Section 4.5.9]{Petrov}, \cite[Section 
 7.2]{Boal}). Moreover, at least in case of some monolayers (see 
 \cite{TKS-Gauss, SiegKoz-Gauss}), $\elbend$ and $\elGauss$  
 satisfy
\begin{equation}\label{eq:constr}
-1<\frac{\elGauss}{\elbend}<0.
\end{equation}
In this paper we are concerned with the variational approximation of 
$\Helf$, under condition \eqref{eq:constr} and with $H_0=0$;
in Section \ref{bombemerda} we briefly discuss how to relax these two
constraints.  
In this respect we note that, for any given $H_0 \in \R$, a 
condition ensuring compactness and lower 
semicontinuity 
of $\Helf$ in a reasonable topology 
 (see Theorem 
\ref{the:Helf-lsc} and Remark \ref{rem:BiKaMi})
is
the existence of  two positive numbers $\constensurecompact$ and $\lambda$ such that
\begin{equation*}
\frac{\elbend}{2}\left(H_{\partial E}-H_0\right)^2+\elGauss K_{\partial E}\geq
 \constensurecompact
\left\vert \sff_{\partial E}\right\vert^2-\lambda,
\end{equation*}
where $\sff_{\partial E}$ denotes the second fundamental form of 
$\partial E$. Such a condition is equivalent to the constraint 
$-2<\elGauss/\elbend<0$ (see Section \ref{Del-puzzone}), which is 
trivially satisfied when \eqref{eq:constr} holds. 

Recently several  authors have used  
diffuse interfaces approximations in order to develop
 efficient numerical simulations for a 
number of models involving  $\Helf$  (e.g. see \cite{BiKaMi, DuSpontCurv, DuWill, DuUno, DuDue, Dutopo, 
DuCa, CampHern, CampHern2, DuStokes, Dumulti, FarGar}). 
Analytical results 
 have been carried on, mainly 
 by means of formal asymptotics,  in \cite{DuDue, DuSpontCurv, DuWill, Wango}.    Most of the papers cited above 
concentrate on the approximation of the 
term $\mathpzc H$ which 
(up to minor modifications)
takes the form
\begin{equation}
\mathpzc H_\eps(u):=\frac{1}{\eps}\int_\Omega\left(
\eps\Delta u- \frac{W^\prime(u)}{\eps}-\eps\vert\nabla u\vert 
H_0\right)^2\, dx,
\label{eq:def-Heps}
\end{equation}
where 
$\eps>0$ is a small parameter related to the width of the diffuse interface, and $W\in 
C^2(\R)$ is a double-well potential with two equal
 minima (from now on, 
throughout the paper, we will make the 
choice $W(s):=(1-s^2)^2/4$). Actually, in the case $H_0=0$, it was firstly conjectured in \cite{DG} that functionals
similar to \eqref{eq:def-Heps}  $\Gamma$-converge 
to $\sigma \mathcal H$
as $\eps\to 0^+$, where $\sigma$ is a suitable
positive constant.

 At least in the case $H_0=0$, 
 the choice of the sequence in \eqref{eq:def-Heps}
can be heuristically motivated with the fact that $\mathcal H_\eps$ represents a kind of (rescaled) squared ``
$L^2$-gradient'' of the functional  $\mathpzc 
P_\eps$ defined as
\begin{equation*}
\mathpzc P_\eps(u):=\int_\Omega
\left(
\frac{\eps}{2}\vert \nabla 
u\vert^2+\frac{W(u)}{\eps}
\right)\,dx, \qquad \text{if }u\in H^1(\Omega),
\end{equation*}
and $\mathpzc P_\eps(u):=+\infty$ elsewhere in $L^1(\Omega)$. This, 
together with the  well known results that $\mathpzc 
P_\eps$  
approximate the perimeter 
functional as $\eps \to 0^+$ (see \cite{MM,Braides-book}), 
and that the ``$L^2$-gradient'' of the 
perimeter is 
formally
given by the mean curvature operator, furnishes a (very) heuristic justification for the choice of $\mathpzc H_\eps$.

The aim of this paper is twofold: we want to 
propose a diffuse interface approximation of 
$\mathpzc K$ which slightly differs from those proposed
until now 
(see \cite{Dutopo, DuCa} and Remark \ref{rem:Du-2}). 
Moreover, we want to prove a rigorous
convergence result  for our approximating sequence within the framework of $\Gamma$-convergence, under the assumptions that $H_0=0$, and 
provided the parameters $\elbend,\,\elGauss$ satisfy \eqref{eq:constr}.
  
In order to define the approximating 
functionals we need some notation.
For every $u\in C^2(\Omega)$ we define the vector field $\nu_u\in 
L^\infty(\Omega)$ by $\nu_u:=\nabla u/\vert\nabla u\vert$ whenever 
$\nabla u\neq 0$ and $\nu_u:=\mathbf e$ on $\{\grad u=0\}$, where 
$\mathbf{e}$ 
is an arbitrary unit vector (to fix the notation from now on we will choose $\mathbf e=\mathbf e_3$, $\mathbf e_3$ being the third element of the canonical basis of $\R^3$).
Then, denoting by $\vert\cdot\vert$ the norm of a matrix as defined in \eqref{amantedellosplit}, 
we propose 
to approximate  $\mathpzc K$ 
with the functionals
$\mathpzc K_\eps$ 
defined as
\begin{align}
\mathpzc K_\eps(u):=&
\frac{1}{2\eps}\int_\Omega
\left[
\left(\eps\Delta u-\frac{W^\prime(u)}{\eps}\right)^2-\left\vert \eps\nabla^2 u-\frac{W^\prime(u)}{\eps}\nu_u\otimes\nu_u \right\vert^2
\right]
\,dx\notag
\\
=&\frac{1}{\eps}\int_\Omega\sum_{1\leq i<j\leq 3}\det\left[\eps\nabla^2 u-\frac{W^\prime(u)}{\eps}\nu_u\otimes\nu_u\right]_{ij}\,dx,
\label{eq:def-Keps}
\end{align}
when $u\in C^2(\Omega)$ and $+\infty$ elsewhere in $L^1(\Omega)$,
where, for a $3\times 3$-matrix $M$, $M_{ij}$ stands
for its $ij$-th principal minor.
Eventually, as an approximation of $\Helf$, 
if $\mathpzc H_\eps$ is as in \eqref{eq:def-Heps}
with $H_0=0$,
we consider
\begin{gather}\label{veroe'benpindemonte}
\mathpzc W_\eps(u)
:=\frac{\elbend}{2}\mathpzc H_\eps(u)+\elGauss\mathpzc K_\eps(u)
\\
=\int_\Omega
\left\{
\frac{\elbend+\elGauss}{2\eps}\left[\mathrm{tr}\left(\eps\nabla^2 u-\frac{W^\prime(u)}{\eps}\nu_u\otimes\nu_u 
\right)
\right]^2
-\frac{\elGauss}{2\eps}\left\vert \eps\nabla^2 u-\frac{W^\prime(u)}{\eps}
\nu_u\otimes\nu_u \right\vert^2\right\}\,dx.
\notag
\end{gather}
We can roughly summarize our main results as follows. 
Suppose that  
\eqref{eq:constr} holds, that $H_0=0$,
 and let $\{u_\eps\}_{\eps}\subset C^2(\Omega)$ satisfy
\begin{equation}\label{eq:ener-bound}
\sup_{0<\eps<1}\mathpzc P_\eps(u_\eps)<+\infty,\qquad 
\sup_{0<\eps<1} \mathpzc W_\eps(u_\eps)< +\infty.
\end{equation}
Then 

\medskip
\noindent(\textit{Compactness},  see 
Theorems \ref{the:Gamma-liminf-cond} and \ref{theo:se2mni}). 
Up to a (not relabelled) subsequence, there exists
a function $u=2\chi_E-1\in BV(\Omega,\{-1,1\})$ such that $\lim_{\eps\to 
0^+}u_\eps=u$ in $L^1(\Omega)$. Furthermore, the measures $\mu_{u_\eps}$ 
associated with the density of the functionals $\mathpzc P_\eps(u_\eps)$ (see 
\eqref{eq:def-mu-eps})  concentrate, as $\eps\to 0^+$, on a generalized 
surface $
\rectifiableset
\supseteq \Omega \cap \partial E$,
 for which a 
weak notion of second 
fundamental form is defined.
Actually, for almost every $s\in (-1,1)$ the oriented varifolds 
associated with the level sets $\{u_\eps=s\}$ converge to the same limit.

\medskip

\noindent (\textit{Lower bound}, see Theorem \ref{the:Gamma-liminf-cond}). 
The $\liminf_{\eps\to 0^+}\mathpzc W_\eps(u_\eps)$ is bounded 
from below by a suitable 
positive constant $\surftens$ times the value of (a suitable extension of) $\Helf$ 
evaluated on $\rectifiableset$. In particular if $E$ has 
$C^2$-boundary in $\Omega$ we have
\begin{equation}\label{brutt}
\liminf_{\eps\to 0^+} \mathpzc W_\eps(u_\eps)\geq\surftens\Helf(E).
\end{equation}

\medskip

\noindent(\textit{Upper bound}, see Theorem \ref{the:Gamma-limsup}).  For 
every bounded open set $E\subset \Omega$  with $C^2$-boundary  
there exists a sequence $\{u_\eps\}_{\eps}\subset C^2(\Omega)$ such 
that $\lim_{\eps\to 0^+}u_\eps\to 2\chi_E-1$ in $L^1(\Omega)$, and 
$\lim_{\eps\to 0^+}\mathpzc W_\eps(u_\eps)=\surftens\Helf(E)$. 
 
\medskip

\noindent($\Gamma(L^1)$-\textit{Limit on smooth points}, see
Corollary \ref{cor:Gamma-conv-smooth-points}). By 
the $L^1(\Omega)$-lower semicontinuity of $\Helf$ (see
Theorem \ref{the:Helf-lsc})
we can conclude that if the bounded set $E$ 
has $C^2$-boundary in $\Omega$, then 
$$
\Gamma(L^1)-\lim_{\eps\to 0^+} \mathpzc W_\eps(u)=\surftens\Helf 
(E). 
$$

As we already said,
in \cite{Dutopo, DuCa} slightly different approximations of the Gaussian curvature have been proposed and used in numerical experiments to retrieve topological informations for the diffuse interface. 
The functional $\K_\eps$ 
in \eqref{eq:def-Keps}
might have some advantages, at least 
from the analytical point of view. 
Firstly
 $\mathpzc W_\eps$ can be expressed in terms of the trace and the norm of 
$\eps\nabla^2 u- \frac{W^\prime(u)}{\eps} \nu_u\otimes\nu_u$, and
 for every $x_0\in\Omega$ such that $\nabla u(x_0)\neq 0$, the 
matrix $\eps\nabla^2 u(x_0)- \frac{W^\prime(u(x_0))}{\eps}\nu_u(x_0)
\otimes\nu_u(x_0)$ 
has an explicit relation with 
the second fundamental form of the level line  $\{u=u(x_0)\}$ times $\vert\nabla u(x_0)\vert$ (see \eqref{eq:acculo?}). Secondly, if \eqref{eq:constr} is satisfied, from \eqref{eq:ener-bound} we can derive the 
bound
\begin{equation*}
\sup_{0<\eps<1}\frac{1}{\eps}\int_\Omega\left\vert\eps\nabla^2u_\eps-\frac{W^\prime(u_\eps)}{\eps}
\nu_{u_\eps}\otimes\nu_{u_\eps}\right\vert^2\,dx<+\infty.
\end{equation*}
{}From this latter relation we can deduce two rather interesting further 
properties. The first is  that, as already stated above, the energy 
measures $\mu_{u_\eps}$  concentrate on a generalized surface 
with second 
fundamental form in $L^2$ (namely a \textit{Hutchinson's 
curvature varifold}, see Lemmata \ref{lem:Kazumba} 
and \ref{lem:conv-to-Hutch}).
As a consequence  we get better regularity for the 
limit of the $\mu_{u_\eps}$ with respect to the case when only a 
uniform bound on 
$\mathpzc H_\eps(u_\eps)$ is available;
indeed, 
 under this latter uniform bound,  
the measures $\mu_{u_\eps}$ concentrate on a \textit{rectifiable 
integral Allard's varifold} with squared integrable generalized mean curvature
(see \cite{RS, ToneYuko}, and
 Appendix \ref{luponeassassino} for 
the definitions of varifold and curvature varifold). The second property is an 
improved convergence to zero 
of the 
discrepancies $\xi^\eps_{u_\eps}$ defined in \eqref{eq:def-xi-eps}. 
In fact, we obtain that 
$\lim_{\eps\to 0^+}\| \frac{\eps}{2} \vert\nabla u_\eps\vert^2-
\frac{W(u_\eps)}{\eps}\|_{L^p(\Omega)}=0$,  for every $p\in[1,3/2)$ (see Proposition \ref{lem:conv-discr-BV}). Let us 
stress that the 
improved convergence of the discrepancies may indicate a good behaviour of $\mathpzc W_\eps$ in 
numerical experiments. 
Indeed, given $\{u_\eps\}_{\eps}\subset C^2(\Omega)$ such that $\lim_{\eps\to 0^+}u_\eps= 2\chi_E-1$ in $L^1(\Omega)$, the 
condition $\frac{\eps}{2}\vert\nabla u_\eps\vert^2-
\frac{W(u_\eps)}{\eps}= O(\eps)$ is one
of the characteristics for a sequence to be a 
``good'' recovery sequence (like, for example, 
the one constructed in Theorem \ref{the:Gamma-limsup}). In other words,
one of the properties that suggests  a ``good'' convergence to the sharp 
interface functional is that  $\frac{\eps}{2}\vert\nabla 
u_\eps\vert^2- \frac{W(u_\eps)}{\eps}$ 
vanishes rapidly enough as $\eps\to 0^+$. 
In numerical applications, a penalizing term of the form 
$\|\frac{\eps}{2}
\vert\nabla u_\eps\vert^2- \frac{W(u_\eps)}{\eps}\|^p_{L^p(\Omega)}$ is 
often added to the diffuse interface functional 
to force such a ``fast'' decay of $\vert\xi^\eps_{u_\eps}\vert$.

Let us conclude by 
remarking the fact that, although an approximation via diffuse 
interfaces seems to be reasonable for numerical purposes, our result does 
not establish any physical derivation of the Helfrich's energy as a 
mesoscale limit, as for example it has been recently done in \cite{PellRoeg}.

The paper is organized as follows. In Section \ref{smerd} we fix some notation, recall some basic definitions from differential geometry and briefly comment 
on the definition of $\mathpzc W_\eps$, 
as well as 
on the relation of $\mathpzc K_\eps$ with  \cite{Dutopo,DuCa}. In Section \ref{puppappera} we summarize  the main results proved in \cite{RS}, that represent one of the pillars on which 
our paper rests. In Section \ref{sec:stat-main}
we state 
our main results. The proofs are postponed to Sections \ref{sec:cavallona}-\ref{WUgoCenci}. In Section \ref{bombemerda} we collect some 
additional results, and
 we show how the assumptions on the parameters $\elbend,\,\elGauss,$  can be weakened; 
we briefly discuss the possibility of proving a full $\Gamma$-convergence result and the problems arising in the case $H_0\neq 0$. 
Eventually in Appendices \ref{trappolapernani}-\ref{luponeassassino} we collect some definitions and results on measure-function pairs  and
geometric measure theory, needed in the proofs of the main results.

\section{Notation}\label{smerd}
\subsection{Linear algebra}\label{sublin}
We endow the space of  the $(3\times 3)$ 
matrices 
$M=(m_{ij})\in\R^{3\times 3}$
(resp. $3^3$ tensors $T=(t_{ijk})\in\R^{3^3}$) with 
the norm
\begin{equation}\label{amantedellosplit}
\vert M\vert^2 := {\rm tr}(M^T M) =\sum_{i,j=1}^3 (m_{ij})^2 \qquad 
\left({\rm 
resp.} 
~\vert T\vert^2:=\sum_{i,j,k=1}^3(t_{ijk})^2 \right),
\end{equation}
where $M^T$ is the transposed of $M$. 

If $M\in\R^{3\times 3}$ is symmetric, 
$O= (o_{il})\in O(3)$  and
$D = {\rm diag}(d_{11}, d_{22}, d_{33})$ are such that 
$M=O^TDO$, then
$\vert M\vert^2= \mathrm{tr}(O^TD^2O)
= \sum_{l=1}^3 
(d_{ll})^2
\sum_{i=1}^3 
(o_{li})^2 
=\sum_{l=1}^3
(d_{ll})^2$.
Moreover, still for a symmetric matrix $M \in \R^{3\times 3}$, we have 
 $\frac{1}{2} \left[ (\mathrm{tr}(M))^2 - \mathrm{tr}(M^TM) \right] =
 \sum_{1\leq i<j\leq 3} \det(M_{ij})$, where $M_{ij}$ is the 
$ij$-principal minor of $M$.

\begin{remark}\rm 
If $P\in\R^{3\times 3}$ is a (symmetric) orthogonal projection matrix onto some 
subspace of $\R^3$ and $M$ is 
symmetric, then
\begin{equation}\label{eq:Kill-Bill}
\vert P^T M P\vert^2\leq\vert M\vert^2.
\end{equation}
Indeed
\begin{equation}\label{mer}
\vert 
P^TMP\vert^2=
\sum_{j=1}^3 
\sum_{i=1}^3\left(\sum_{l=1}^3 p_{il}\left(\sum_{k=1}^3 m_{lk} 
p_{kj}\right) 
\right)^2 
=\sum_{j=1}^3\left\vert P\left(MP\right)^{(j)}\right\vert^2,
\end{equation}
where the column vector $\left(MP\right)^{(j)}\in \R^3$ has components
$(\sum_{k=1}^3 m_{1k}p_{kj}, \sum_{k=1}^3 m_{2k}p_{kj},
\sum_{k=1}^3 m_{3k} p_{kj})$, and $\vert \cdot \vert$ 
on the right hand side of 
\eqref{mer} is the  euclidean norm of a vector. 
Since $P$ is an orthogonal projection we have
\begin{equation*}
\vert P^TMP\vert^2
\leq  \sum_{j=1}^3\left\vert \left(MP\right)^{(j)}\right\vert^2
=\sum_{i,j=1}^3\left(\sum_{k=1}^3 m_{ik} p_{kj}\right)^2
=\sum_{i=1}^3\left\vert\left(M\right)_{(i)}P\right\vert^2,
\end{equation*}
where 
$\left(M\right)_{(i)}=(m_{i1},m_{i2},m_{i2})\in\R^3$.
Using again the fact that $P$ is a projection we have
\begin{equation*}
\vert P^TMP\vert^2\leq 
\sum_{i=1}^3\left(\left(M\right)_{(i)}\right)^2
=\sum_{i=1}^3\sum_{j=1}^3 \left(
m_{ij}\right)^2=\vert M\vert^2.
\end{equation*}
\end{remark}

By $G_{2,3}$ (resp.  $G^0_{2,3}$) we denote the Grassmannian 
of the unoriented
$2$-planes in $\R^3$ (resp. the Grassmannian 
of the oriented $2$-planes in $\R^3$). 

We denote by $\q$ the standard $2$-fold covering map $\q:G^0_{2,3}\to G_{2,3}$. 
We often identify $G^0_{2,3}$ with 
the set of simple unit $2$-vectors  $\tau\in\Lambda_2(\R^3)$.  
Moreover
$$
\star:\Lambda^1(\R^3)\to \Lambda_2(\R^3)
$$
denotes the 
 Hodge operator. Often vectors and covectors will be identified. 
For every $\tau\in G^0_{2,3}$ 
we define $\nu^\tau\in\R^3\simeq\Lambda^1(\R^3)$ as 
the unique unit vector such that $\star\nu^\tau=\tau$.
 
We endow $G_{2,3}$ with the distance induced by the norm 
$\vert S\vert$,
where $S$ is the matrix associated with the orthogonal 
 projection of $\R^3$ onto $S\in 
G_{2,3}$. Moreover, for every open set $\Omega\subseteq\R^3$ we let 
$G_2(\Omega):=\Omega\times G_{2,3}$, endowed with the product 
distance.

In the same way, we endow $G^0_{2,3}$ with the distance induced by 
$\vert\tau\vert$,
 where $\tau$ is the simple unit 
$2$-vector associated with $\tau\in G^0_{2,3}$. 
Moreover, for every open set
$\Omega\subseteq\R^3$ we let $G^0_2(\Omega):=\Omega\times G^0_{2,3}$, 
endowed 
with the product distance.
Finally, we let 
$\mathbb S^2:= \{\xi\in\R^3:~\vert\xi\vert=1\}$, and we denote
by $\triangle$ the symmetric difference between sets.

\subsection{Differential Geometry}\label{prel:diff-geom}
Let $\Sigma$ be a smooth, compact oriented surface without boundary 
embedded in $\R^3$. 
If $x\in\Sigma$, we denote by $P_\Sigma(x)$ the orthogonal
projection onto the tangent plane $T_x\Sigma$ to $\Sigma$ at $x$. 
Often we identify the linear operator $P_\Sigma(x)$  with the symmetric 
$(3\times 3)$-matrix ${\rm Id}-\nu_x\otimes\nu_x$ where $x \to \nu_x \in 
(T_x\Sigma)^\perp$ is a smooth unit 
covector field 
orthogonal to $T_x\Sigma$.

Let us recall that, when $\Sigma$ is given as a level surface $\{v=t\}$ 
of a smooth function $v$ such that $\nabla v\neq 0$ on $\{v=t\}$,  we 
can take at $x \in \{v=t\}$
\begin{equation*}
\nu_x=\frac{\nabla v(x)}{\vert\nabla v(x)\vert},\qquad P_\Sigma(x)
={\rm Id}-\frac{\nabla v(x)\otimes\nabla v(x)}{\vert\nabla v(x)\vert^2}.
\end{equation*}
The second fundamental form ${\bf B}_\Sigma$ of $\Sigma$ has the 
expression
\begin{equation*}
\sff_{\Sigma}=
\Big( P_\Sigma^T 
\frac{\nabla^2v}{\vert\nabla v\vert}
P_\Sigma\Big)
\otimes \frac{\nabla v}{\vert\nabla v\vert},
\end{equation*}
where $P_\Sigma^T = (P_\Sigma)^T$.
The definition of $\sff_\Sigma$ depends only on $\Sigma$ and not 
on the particular choice of the function $v$. Moreover 
$\sff_\Sigma(x)$,
if restricted to $T_x \Sigma$ and 
considered as a bilinear map from $T_x \Sigma \times T_x \Sigma$ 
with values in $(T_x \Sigma)^\perp$, 
 coincides with the usual notion of second fundamental form.
By 
$$
\Hv_\Sigma(x)=(H_1(x),H_2(x),H_3(x)) = {\rm tr}\Big( P_\Sigma^T 
\frac{\nabla^2v}{\vert\nabla v\vert}
P_\Sigma\Big) \nu_x,
$$
we denote the mean curvature vector of $\Sigma$ at $x \in \Sigma$. 
We define
the (scalar) mean curvature of $\Sigma$ at $x$ with respect to $\nu_x$ as
$$
H_{\Sigma}(x):=\Hv_\Sigma(x)\cdot\nu_x.
$$ 
Notice that $\Hv_\Sigma$ does not depend on the choice of 
$\nu$, while the sign of $H_{\Sigma}$ does. Observe also that 
$H_{\Sigma}$ is the sum of the two principal curvatures of $\Sigma$: 
sometimes $H_{\Sigma}$ is also referred to as the total curvature. 
When  $\Sigma=\partial E$,  where $E\subset\R^3$ is open and bounded, 
we define $\nu_{\partial E}$ to be the interior normal to $\partial 
E=\Sigma$ and   $H_{\partial E}
:=\Hv_{\partial E}\cdot\nu_{\partial E}$, which turns out to 
 positive on convex surfaces.

Let us also define $A^\Sigma(x):=(A^\Sigma_{ijk}(x))_{1\leq i,j,k\leq 
3}\in\R^{3^3}$ as
\begin{equation}\label{eq:def-Aijk-Sigma}
A^\Sigma_{ijk}=\delta^\Sigma_i P_{\Sigma jk} \qquad {\rm on} ~ \Sigma,
\end{equation}
where $\delta_i^\Sigma := P_{\Sigma ij} \frac{\partial}{\partial x_j}$.

To better understand definition \eqref{eq:def-Aijk-Sigma}, 
it is useful to recall the 
links netween 
$\sff_\Sigma$ and $A^\Sigma$ (see \cite[Proposition 2.3]{Mant}).

\begin{proposition}\label{prop:A-vs-B}
Set $A = A^\Sigma$, ${\bf B} = {\bf B}_\Sigma$
and ${\bf H} = {\bf H}_\Sigma$.
For $i,j,k \in \{1,2,3\}$
the following relations 
hold:
\begin{align}
& B^k_{ij}=P_{jl}A_{ikl},
\label{Bjork}
\\
& A_{ijk}=B^k_{ij}+B^j_{ik},
\label{eq:AvsB}
\\
& \Hv_i=A_{jij}=B^j_{ji}+B^i_{jj}.
\end{align}
\end{proposition}
The next proposition shows some of the relations between 
the curvatures of $\Sigma$ 
and the derivatives of the signed distance function from $\Sigma$ itself.
\begin{proposition}\label{prop:dist-VS-curv}
Let $E$ be a bounded open subset of $\R^3$ with $C^2$-boundary.
Then there exists an open neighborhood $U$ of $\partial E$  such that, 
denoting by $d:U\to \R$ the signed distance from $\partial E$ positive inside $E$, we have
$d\in C^2(U)$ and, for $y \in U$ and 
 $\pi(y):=y-d(y)\nabla d(y)\in\partial E$ the unique orthogonal projection
point of $y$ on $\partial E$, 
\begin{gather}\label{babbo}
 \Delta d(y) = H_{\partial E}(\pi(y))+o(d(y))
\\
 \sum_{1\leq i<j\leq 3}\det\left([\nabla^2 d(y)]_{ij}\right)= 
K_{\partial E}(\pi(y))+ o(d(y)),
\label{eq:Gauss-VS-dist}  
\end{gather}
where $o(t)\to 0$ as $t\to 0$.
\end{proposition}
\proof
It is 
well known (see for example \cite{G}) 
that $d$ is of class $C^2$ in a suitable 
tubular neighborhood $U$ of $\partial E$ where
$\pi$ is single valued, and 
moreover that, for every $y\in U$, the eigenvalues of $\nabla^2 d(y)$ are 
$$
\lambda_1(y) =\frac{k_1(\pi(y))}{1-
d(y)k_1(\pi(y))},\qquad\lambda_2(y) =\frac{k_2(\pi(y))}{1-
d(y)k_2(\pi(y))},\qquad \lambda_3(y)=0,
$$
where
 $k_1(x),k_2(x)$ are the principal curvatures of $\partial E$ at $x$. 
Then \eqref{babbo} follows, and
\begin{equation*}
\begin{aligned}
 \sum_{1\leq i<j\leq 3}\det\left([
\nabla^2 d(y)]_{ij}\right)
= &
\frac{\Big({\rm tr}(\nabla^2 d(y))\Big)^2
-\vert\nabla^2 d(y)\vert^2}{2}
\\
= & \lambda_1(y)\lambda_2(y)
=  K_{\partial E}(\pi(y))+o(d(y)).
\end{aligned}
\end{equation*}
\endproof
\subsection{The Helfrich's Functional $\Helf$.}
Throughout the paper  $\Omega\subseteq\R^3$ is an open connected
set with smooth boundary ($\Omega = \R^3$ is allowed). If $E \subseteq 
\R^3$, $\chi_E$ is the characteristic function of $E$ equal to 
$1$ on $E$ and $0$ elsewhere.
Let $E\subseteq \Omega$ be an open set. We say that $E$ has 
$C^k$-boundary in $\Omega$ ($k\in\N\cup\{\infty\}$) if for every 
$x\in \Omega \cap \partial E$ the set $\Omega \cap \partial E$ 
can be  written, locally around $x$, as the graph of a $C^k$ function,
and $\Omega \cap E$ is locally the subgraph of the 
same function. 

By assumption \eqref{eq:constr} it follows that
$\frac{\elbend-\elGauss}{\elbend+\elGauss}$ is a 
positive
real number. We set
\begin{equation}\label{eq:def-l2}
l^2:=(H_0)^2\elbend ~
\frac{\elbend-\elGauss}{2(\elbend+\elGauss)}.
\end{equation}

\noindent We claim that, 
whenever $E$ is bounded with smooth boundary in $\Omega$, then
\begin{equation*}
\Helf(E)\geq\,-l^2\Ha^2(\Omega \cap \partial E).
\end{equation*}
To prove the claim, write 
\begin{equation*}
\begin{split}
& \Helf(E)
\\
=& \int_{\Omega \cap \partial E}
\left[-\frac{\elGauss}{2}\vert \sff_{\partial 
E}\vert^2
+\left(\frac{\elbend + \elGauss}{2}\right) (H_{\partial E})^2+
\elbend H_0 H_{\partial E}+\frac{\elbend}{2} (H_0)^2\right]\,d\Ha^2.
\end{split}
\end{equation*}
If  $\alpha 
:= \frac{\elbend+\elGauss}{2}
>0$, $\beta := \elbend H_0$, and $\gamma  := \elbend \frac{(H_0)^2}{2}$,
since
$\alpha t^2 + 
\beta t + 
\gamma \geq \frac{\alpha}{2} t^2 - l^2$ for any $t \in \R$ 
and
$l^2= \frac{\beta^2}{2\alpha} - \gamma$, we have
the inequality
\begin{equation}\label{eq:Helf-bded-below}
\Helf(E)
\geq \int_{\Omega \cap \partial E}\left[-\frac{\elGauss}{2}\vert \sff_{\partial 
E}\vert^2
+\frac{(\elbend+\elGauss)}{4} (H_{\partial E})^2-l^2\right]\,d\Ha^2.
\end{equation}
Thanks to \eqref{eq:constr}, the first
two addenda inside the integral on the right hand side of 
\eqref{eq:Helf-bded-below}
are nonnegative, hence the claim follows. 

\subsection{Definitions of $\muepsu$, $\tildemuepsu$, $\xiepsu$, 
$R_u^\eps$,
$\BBu$,
$A^u$, 
$\Vepsu$, 
$\Vzeroepsu$,
$\fepsu$, 
$\calBepsu$ and
$\Hepsu$.
}
We set
\begin{equation*}
W(r):=\frac{1}{4} (1-r^2)^{2}, \qquad r \in \R,
\end{equation*}
and
\begin{equation}\label{dito}
\surftens:=\int_{-1}^1\sqrt{2W(s)}\,ds.
\end{equation}
If $\gamma(s):=\tanh(s)$ we have
$\ddot \gamma=\frac{d}{ds}
(W(\gamma))$,
\\
\begin{gather*}
\int_\R\vert\dot\gamma\vert^2\,ds=\int_\R 2W(\gamma)\,ds
=\surftens,
\end{gather*}
and
\begin{equation}\label{eq:def-surf-tens}
\surftens=\min\left\{\int_\R\left(\frac{\vert \dot v\vert^2}{2}+
W(v)\right)\,ds:\, v\in H^1_{{\rm loc}}(\R),\,
\lim_{s\to\pm\infty}v(s)=\pm1\right\}.
\end{equation}

For $u \in \C^2(\Omega)$ and 
$\LL^3$ 
the Lebesgue measure in $\R^3$,
we define the following Radon measures:
\begin{align}
& \muepsu:=\left(\frac{\eps}{2}\vert\nabla 
u\vert^2+\frac{W(u)}{\eps}\right) \LL^3 \res\Omega,
\label{eq:def-mu-eps}
\\
& \tildemuepsu:=\eps\vert\nabla 
u \vert^2\LL^3 \res\Omega,
\label{eq:def-mu-tilde-eps}
\\
& \xiepsu:=\left(\frac{\eps}{2}\vert\nabla 
u\vert^2-\frac{W(u)}{\eps}\right)\LL^3 \res\Omega,
\label{eq:def-xi-eps}
\end{align}
where $\res$ is the restriction. $\xiepsu$ is usually
called discrepancy measure, while $\muepsu$ is 
the density of the Allen-Cahn functional $\mathcal P_\eps$.
With a small abuse of notation, when necessary we still
denote by $\xiepsu$ the density of the discrepancy measure, i.e.,
$\xiepsu = \frac{\eps}{2} \vert \grad u\vert^2 - \frac{W(u)}{\eps}$.
Note that 
\begin{equation}\label{grad:discr}
\grad \xiepsu = \eps \grad^2 u \grad u - \frac{W'(u)}{\eps} 
\grad u.
\end{equation}

For $u \in C^2(\Omega)$ define $R^\eps_{u}:G_2(\Omega)\to\R^3$ as
\begin{equation}\label{erbarmedich}
R^\eps_{u}(x,S)=R^\eps_{u}(x):=
\frac{1}{\eps\vert\nabla u(x)\vert^2}
\nabla \xi^\eps_{u}(x),
\end{equation}
with the convention that $R^\eps_{u}:= 0$ 
on the set $\{\nabla u=0\}$.

Let $u \in C^2(\Omega)$.  We will often look at geometric 
properties
of the \textit{ensemble of the level sets} of $u$. 
We 
define
\begin{equation}\label{eq:def-proj-nablaperp}
\niu:=\frac{\nabla u}{\vert\nabla u\vert},\qquad 
\proju:={\rm Id}-\niu \otimes\niu, \qquad
P^u_{ij} = \delta_{ij} - (\nu_u)_i (\nu_u)_j,
\end{equation}
on $\{\grad u\neq 0\}$ and  $\niu:=\mathbf{e}_3,\, \proju:={\rm 
Id}-\mathbf{e}_3\otimes\mathbf{e}_3$ on $\{\grad u=0\}$. Moreover we 
define 
the second fundamental form of the ensemble of the level sets of $u$ 
by
\begin{equation}\label{eq:def-sff-lev-surf}
\sff_u=\frac{(\proju)^T\nabla^2 u \proju}{\vert\nabla u\vert}\otimes\niu,
\end{equation}
on $\{\grad u\neq 0\}$ 
and $\sff_u:=\otimes^3\mathbf{e}_3$ on $\{\grad u=0\}$. 
Similarly we define
\begin{equation}\label{def:Aijk}
A_{ijk}^u:=-\proju_{il}\big[\partial_{l}((\niu)_j(\niu)_k)\big],
\end{equation}
on $\{\nabla u\neq 0\}$ and 
$A^u:=\otimes^3\mathbf{e}_3$ on $\{\nabla u=0\}$.

It will be convenient to consider $\sffu$ and $A^u$ 
as  defined on 
$G_2(\Omega)$ (resp. on $G^0_2(\Omega)$) by 
$\sff_u(x,S):=\sff_u(x)$,
$A^u(x,S):=A^u(x)$
 (resp. $\sff_u(x,\tau):=\sff_u(x)$,
 $A^{u}(x,\tau):=A^{u}(x)$
).

By $V_{u}$ (resp. $V^{0}_u$) we 
denote the varifold (resp. oriented varifold)
\begin{align}
\Vepsu(\phi)=\surftens^{-1}\int\phi(x,\proju)~d\tildemuepsu \qquad \forall 
\phi\in C^0_c(G_2(\Omega)),
\label{ehhh???}
\\
\Vzeroepsu
(\phi)=\surftens^{-1}\int\phi(x,\star \niu)~d\tildemuepsu\qquad 
\forall \phi \in C^0_c(G^0_2(\Omega)),
\label{giocare}
\end{align}
see Appendix \ref{luponeassassino}.

We also set 
\begin{equation}\label{eq:def-chem-pot}
\fepsu:=\eps\Delta u-\frac{W^\prime(u)}{\eps}.
\end{equation}
\begin{definition}
Let $u \in \C^2(\Omega)$ and $x\in\Omega$. We  
define 
\begin{equation}\label{eq:def-Aueps}
\calBepsu(x):=
\begin{cases} \displaystyle
\frac{1}{\eps\vert\nabla u(x)\vert} 
\left(\eps\nabla^2u(x)-\frac{W^\prime(u(x))}{\eps}
\niu(x)\otimes\niu(x)\right) & {\rm if}~
\nabla u(x)\neq 0,
\\
0 & {\rm otherwise},
\end{cases}
\end{equation}
\begin{equation}\label{eq:def-Hueps}
\Hepsu(x):=\mathrm{tr}\Big(\calBepsu(x)\Big)=\begin{cases}
\displaystyle \frac{\fepsu(x)}{\eps\vert\nabla 
u(x)\vert} & {\rm if ~}\nabla u(x)\neq 0,
\\
0 & {\rm otherwise}.
\end{cases}
\end{equation}
\end{definition}
We can informally think of 
$\calBepsu\otimes\niu$ and $\Hepsu\niu$ 
as the 
{\it approximate} second 
fundamental 
form and the approximate  mean curvature 
vector of the level sets of $u$, respectively.

Note  that 
$$
R^\eps_{u}
=\calBepsu\frac{\nabla u}{\vert\nabla u\vert} \qquad {\rm on}~
\{\nabla u\neq 0\}.
$$
%

\subsection{The functionals $\mathpzc W_\eps$}
We recall that 
our approximating sequences of functionals is defined in 
\eqref{veroe'benpindemonte},  
  where $\mathpzc H_\eps, \mathpzc K_\eps$ are as in \eqref{eq:def-Heps}, 
\eqref{eq:def-Keps}. 

Observe that 
\begin{equation*}
\int (\Hepsu)^2\,d\tildemuepsu
\leq\mathpzc 
H_\eps(u),
\end{equation*}
with equality if  $\LL^3(\{\fepsu\neq 0\}\cap
\{\nabla u= 0\})=0$, and
\begin{equation*}
\int_\Omega \vert \calBepsu
\vert^2\,d\tildemuepsu
\leq \frac{1}{\eps}
\int_\Omega\left\vert\eps\nabla^2 
u-\frac{W^\prime(u)}{\eps}\niu\otimes\niu\right\vert^2\,dx
,
\end{equation*}
with equality if
\begin{equation*}
\LL^3\left(\left\{\eps\nabla^2 
u-\frac{W^\prime(u)}{\eps}\niu\otimes\niu\neq 0\right\} 
\cap\left\{\nabla u= 0\right\}\right)=0.
\end{equation*}

\noindent Moreover
\begin{equation*}
\int_\Omega \frac{(\Hepsu)^2-\vert \calBepsu
\vert^2}{2}\,d\tildemuepsu=\int_{\{\nabla u\neq 0\}}
\sum_{1\leq i<j\leq 3}^3\det\left( \left[\calBepsu
\right]_{ij}\right)\,d\tildemuepsu,
\end{equation*}
where $\left[\calBepsu\right]_{ij}$ is the $ij$-th 
principal 
minor of $B_u^\eps$, and
\begin{equation}\label{eq:princ-min-i}
\begin{split}
\det\left(\left[\calBepsu\right]_{ij}\right)=
&\frac{1}{\eps^2\vert\nabla u\vert^2}\Bigg[\left(\eps\partial^2_{ii}u
-\frac{W^\prime(u)}{\eps}\frac{(\partial_i u)^2}{\vert\nabla 
u\vert^2}\right)
\left(\eps\partial^2_{jj}u
-\frac{W^\prime(u)}{\eps}\frac{(\partial_j u)^2}{\vert\nabla 
u\vert^2}\right)
\\
&-\left(\eps\partial^2_{ij}u
-\frac{W^\prime(u)}{\eps}\frac{\partial_i u\partial_j u}{\vert\nabla 
u\vert^2}\right)^2\Bigg].
\end{split}
\end{equation}

\begin{remark}\label{rem:Du-2}\rm
Let us notice that
\begin{equation*}\label{eq:Du-2}
\begin{split}
&\frac{(\Hepsu)^2-\vert\calBepsu\vert^2}{2}=
\frac{(\fepsu)^2-\mathrm{tr}[(\eps\nabla^2u-\frac{1}{\eps}W^\prime(u)
\niu\otimes\niu)^2]}{2\eps^2\vert\nabla u\vert^2}
\\
=& \frac{1}{2\eps^2\vert\nabla u\vert^2}
\left((\fepsu)^2-\mathrm{tr}\Big[\eps^2(\nabla^2u)^2-2W^\prime(u)\nabla^2 u
~\niu\otimes\niu+\frac{(W^\prime(u))^2}{\eps^2}
\niu\otimes\niu\Big]\right)
\\
=&\frac{\eps^2\left\{
(\Delta u)^2-\mathrm{tr}[(\nabla^2 u)^2]\right\}-2W^\prime(u)(\Delta u-\partial^2_{\niu \niu}u)}{2\eps^2\vert\nabla u\vert^2}
\\
=&\frac{1}{2\eps^2\vert\nabla u\vert^2}\Big\{\eps^2
\mathrm{div}(\Delta u\nabla u-\nabla^2u\nabla u)-
2W^\prime(u)\mathrm{tr}\left[({\rm Id}-\niu \otimes\niu)
\nabla^2 u)\right] \Big\},
\end{split}
\end{equation*}
where we used 
\begin{equation*}
\mathrm{div}(\nabla^2 u\nabla u)=\mathrm{tr}[(\nabla^2u)^2]+
\nabla u\cdot\nabla(\Delta u).
\end{equation*}
Suppose that $\Omega\subset\subset\R^3$ is open, and $u\in C^2(\Omega)$ verifies $\nabla u\equiv 0$ on $\Omega\setminus \Omega^\prime$, for some
 $\Omega^\prime\subset\subset\Omega$. By Sard's Lemma we can find a sequence of $t_k\in \R^+$ such that $t_k\to 0$ as $k\to\infty$, and, setting $N_k:=\{\vert\nabla u\vert>t_k\}$ we have
\begin{gather*}
\partial N_k\subseteq\{\vert\nabla u\vert=t_k\}
\text{ is a smooth, embedded surface},
\\
\lim_{k\to\infty}\LL^3\Big(\{\nabla u\neq  0\}\setminus N_k\Big)=0.
\end{gather*}

Thus we have
\begin{equation*}
\begin{split}
\left\vert\int_{\{\nabla u\neq 0\}}
\mathrm{div}(\Delta u\nabla u-\nabla^2 u\nabla u)\,dx\right\vert=&
\lim_{k\to\infty}\left\vert\int_{N^\eps_k}
\mathrm{div}(\Delta u\nabla u-\nabla^2 u\nabla u)\,dx\right\vert
\\
=&\lim_{k\to\infty}\left\vert\int_{\partial N_k}
(\Delta u\nabla u-\nabla^2 u\nabla u)\cdot \nu_{\partial N^\eps_k}\,d\Ha^2\right\vert
\\
\leq& \lim_{k \to\infty}\|u\|_{C^2}\Ha^2(\partial N_k)t_k=0.
\end{split}
\end{equation*}
Hence
\begin{equation*}
\begin{split}
&\int\frac{(\Hepsu)^2-\vert\calBepsu\vert^2}{2}~d\tildemuepsu
\\
=&\frac{1}{2\eps}\int_{\{\nabla u\neq 0\}}
\Big(
\eps^2\mathrm{div}(\Delta u\nabla u-\nabla^2 u\nabla u)
-2W^\prime(u)~
\mathrm{tr}\left[({\rm Id}
-\niu\otimes\niu)\nabla^2 u)\right] \Big)\,dx
\\
=&-\int_{\{\nabla u\neq 0\}}\frac{W^\prime(u)}{\eps}
\mathrm{tr}
\left[\proju\nabla^2 u\right]\,dx.
\end{split}
\end{equation*}
When $u_\eps(x)=\gamma_\eps(d(x))+g_\eps(x)$, where $\gamma_\eps$ is as in Section \ref{prendodimiraquestalberoinmezzoalgiardino} and $g_\eps\in C^2(\Omega)$ is such that $\|g_\eps\|_{C^2(\Omega)}=O(\eps)$, this formula coincides (up to an error of order $O(\eps)$)  with the one proposed in \cite{DuCa} 
in order to approximate $\mathpzc K$.
\end{remark}
%
\section{Preliminary known results}\label{puppappera}
In this section we recall some recent results about a modified  
conjecture of De Giorgi concerning the variational approximation  of the Willmore 
functional (see \cite{DG}).  More 
precisely, the so-called $\Gamma-\limsup$ inequality has been  proved in 
\cite{BePa} in any dimension on smooth boundaries; in \cite{BMRad} the 
$\Gamma-\liminf$ inequality has been 
proved in any dimension, under a rather strong ansatz  on the $u_\eps$ 
(namely $u_\eps=v_\eps(d)$,  where $d$ is the signed distance from the 
boundary of the limit set). An ansatz-free proof of the $\Gamma-\liminf$ 
inequality
has been given in dimension $2$ and $3$ in \cite{RS},  and independently, but only
in two-dimensions, in \cite{ToneYuko} (by means of 
 a different proof which makes use of  
generalized varifolds  introduced in \cite{MosDG}).

The following theorem has been proved in \cite{RS} and is one 
of the key 
ingredients in the proofs of our results. 

\begin{theorem}\label{the:RS} Let $\{u_\eps\}\subset C^2(\Omega)$ be a
sequence such that 
$$ \sup_{0<\eps<1}\left\{
\muepsueps(\Omega)+\frac{1}{\eps}\int_\Omega\left( \eps\Delta
u_\eps-\frac{W^\prime(u_\eps)}{\eps}\right)^2\,dx \right\} < +\infty. 
$$
Then there exists a subsequence (still denoted by $\{u_\eps\}$)
converging to $u=2\chi_E-1$ in $L^1(\Omega)$, where $E$ is a finite
perimeter set. Moreover 
\begin{itemize} 
\item[(A)]
$\mu^\eps_{u_\eps}\rightharpoonup\mu$ as $\eps\to 0^+$ weakly$^*$ in $\Omega$ 
as Radon measures and $\mu$ 
verifies 
\begin{equation*} \mu \geq \surftens \Ha^2 \res\partial E.
\end{equation*} 
In addition 
\begin{equation}\label{eq:the-RS-van-discr}
\lim_{\eps\to 0^+}\vert\xiepsueps\vert=0\qquad {\rm ~as~ Radon~measures},
\end{equation} 
where $\vert \xiepsueps\vert$ 
denotes the total variation of the measure 
$\xiepsueps$, 
and hence 
\begin{equation}\label{eq:lim-tuttuguale}
\mu = \lim_{\eps\to 0^+}\muepsueps=\lim_{\eps\to 
0^+}
\tildemuepsueps
=\lim_{\eps\to 0^+}\frac{2W(u_\eps)}{\eps}\LL^3 \res \Omega\qquad {\rm ~as~
Radon~measures.} 
\end{equation} 
\item[(B)] The sequence $\{\Vepsueps\}$ 
 converges in  the varifolds 
sense to an
integral-rectifiable varifold $V\in\mathbf{IV}_2(\Omega)$ with generalized mean curvature
$\Hv_{V}\in L^2(\mu)$ and such that 
$\mu_V = \surftens^{-1}\mu$. 
\item[(C)] For any $Y\in C^1_c(\Omega;\Rn)$ we have 
\begin{equation}\label{eq:the-RS-conv-first-var} 
\surftens\lim_{\eps\to 0^+}\delta
\Vepsueps(Y) =\lim_{\eps\to 0^+}-\int_{\Omega} \fepsueps\nabla 
u_\eps\cdot Y\,dx=
-\int_{\Omega} \Hv_{V}\cdot Y\,d\mu, 
\end{equation}
and 
\begin{equation}\label{eq:the-RS_lsc-var} 
\surftens\int_{\Omega}\vert
\Hv_{V}\vert^2\,d\mu_V \leq 
\liminf_{\eps\to
0^+}\frac{1}{\eps}\int_{\Omega}\left(\eps\Delta u_\eps
-\frac{W^\prime(u_\eps)}{\eps}\right)^2\,dx.
\end{equation}
\end{itemize} 
\end{theorem}

An important point in order to 
establish the $\Gamma(L^1(\Omega))$-convergence of 
$\mathpzc W_\eps$ 
to $\Helf$ is the lower-semicontinuity of $\Helf$ on smooth sets. 
This is the aim of the following theorem, which is a 
consequence of  \cite[Theorem 5.1]{Del}.
\begin{theorem}\label{the:Helf-lsc}
Let $H_0 \in \R$ and
suppose that \eqref{eq:constr} holds.
Let $E\subset \Omega$ be a bounded open set with smooth boundary in $\Omega$. 
Let $\{E_h\}$
be a sequence  of bounded open subsets of $\Omega$ 
 with smooth boundary in $\Omega$, 
such that 
\begin{align}
&\sup_{h\in\N}\Ha^2(\Omega \cap \partial E_h)<+\infty,
\label{eq:area-bound-sharp}
\\
&\lim_{h\to\infty}\LL^3(\Omega\cap(E_h\triangle E))=0.
\label{eq:L1-conv-sharp}
\end{align}
Then 
\begin{equation}\label{eq:Helf-lsc}
\Helf(E) \leq
\liminf_{h\to\infty}\Helf(E_h).
\end{equation}
\end{theorem}

\begin{remark}\label{rem:Sil-barb}
Theorem \ref{the:Helf-lsc}  holds under the weaker assumption 
$-2<\elGauss/\elbend<0$.
\end{remark}

\begin{remark}\label{rem:BiKaMi}\rm
The bound \eqref{eq:area-bound-sharp} is necessary in order
to gain sufficient compactness on the sequence $\{\partial E_h\}$,
since the bound $\sup_h \Helf(E_h) < +\infty$ alone does not
imply any uniform control on the area of $\partial E_h$. 
This is seen with the following example: $\Omega = \R^3$, $H_0 
= 2$, $E_h$ the union, over $n \in \{1,\dots,h\}$,
 of the balls of radius 
$1$ and centered at 
$(2n,0,0)$, so that $\Helf(E_h) = 4\pi^2 \elGauss h <0$.
\end{remark}

\section{Statements of the main results} 
\label{sec:stat-main}
We can now state our $\Gamma$-convergence results.
\begin{theorem}[Equicoercivity and $\Gamma$-liminf 
inequality]\label{the:Gamma-liminf-cond}
Let $H_0=0$ and  suppose that \eqref{eq:constr} holds. 
Let $\{u_\eps\}\subset C^2(\Omega)$ be a sequence satisfying 
\eqref{eq:ener-bound}. Then 
there exists a (not relabelled) 
subsequence satisfying the theses of 
 Theorem \ref{the:RS}. 
Moreover, the 
 varifold $V$ in Theorem \ref{the:RS}
is a curvature
 varifold with generalized second fundamental form ${\bf B}_V$ in 
$L^2$
(see Definition \ref{def:Hutch-var}), and 
\begin{equation}\label{eq:Aueps-to-AVmu-the}
\lim_{\eps\to 0^+}(\Vepsueps, A^{u_\eps})=(V,A_{V})
\end{equation}
as measure-function pairs on $G_2(\Omega)$ with values in $\R^{3^3}$. 
Eventually 
\begin{equation}\label{eq:Gamma-liminf-var}
\liminf_{\eps\to 0^+}\mathpzc W_{\eps}(u_\eps)
\geq 
\surftens\int 
\left[
\frac{\elbend}{2}\vert\Hv_{V}\vert^2+\frac{\elGauss}{2}(\vert\Hv_V\vert^2-\vert 
\sff_{V}\vert^2)\right]\,dV.
\end{equation}
\end{theorem}

\begin{theorem}[$\Gamma$-limsup inequality]\label{the:Gamma-limsup}
Let $H_0=0$ and $E\subset \Omega$ 
be a bounded open set with boundary of class $C^2$.
Then there exists a sequence $\{u_\eps\}\subset C^{2}(\Omega)$ 
such 
that
\begin{align}
& \lim_{\eps\to 0^+}u_\eps=2\chi_E-1 ~{\rm in~ }L^1(\Omega), 
\label{eq:L1-conv-Gammalimsup}
\\
&\lim_{\eps\to 0^+}\mu^\eps_{u_\eps} = \surftens\Ha^{2} \res \partial E
{~\rm as~ Radon~measures},
\label{eq:Rm-conv-Gammalimsup}
\\
&\lim_{\eps\to 0^+}
\mathpzc W_\eps(u_\eps)
=\surftens\Helf(E).
\label{eq:Gammalimsup}
\end{align}
\end{theorem}

As a consequence of Theorems \ref{the:Gamma-limsup},  
\ref{the:Gamma-liminf-cond} and \ref{the:Helf-lsc} we obtain the following
\begin{corollary}[$\Gamma$-limit on 
smooth sets]\label{cor:Gamma-conv-smooth-points}
Let $H_0=0$ and suppose that 
\eqref{eq:constr} holds. Let $E\subset \Omega$ 
be a bounded open set with boundary of class $C^2$.
Then
 \begin{equation}\label{eq:Gamma-conv-smooth-points}
\left[\Gamma(L^1(\Omega))-\lim_{\eps\to 0^+} \mathpzc 
W_\eps\right](2\chi_E-1)
=\surftens\Helf(E).
\end{equation}
\end{corollary}

Next theorem shows that actually from the hypotheses 
of Theorem \ref{the:Gamma-liminf-cond} we can 
prove a stronger compactness result, 
since 
the oriented varifold 
(see Appendix \ref{luponeassassino})
associated with almost every level line converge to the same limit.
\begin{theorem}[Enhanced compactness]\label{theo:se2mni}
Let $H_0=0$ and suppose that \eqref{eq:constr} holds. 
Let $\{u_\eps\}\subset C^2(\Omega)$ be a sequence 
satisfying
\eqref{eq:ener-bound}. Then there exists a (not relabelled)
subsequence  such that 
\begin{itemize}
\item[(A)] the sequence $\{\Vzeroepsueps\}$ converges in the sense of 
oriented 
varifolds  to an oriented 
varifold $V^0\in\mathbf{IV}^0_2(\Omega)$ such that $\q_\sharp V^0=V$, where $V\in\mathbf{IV}_2(\Omega)$ is as in Theorem \ref{the:Gamma-liminf-cond}.
\item[(B)]  For every $\psi\in C^1_c(\Omega\times\mathbb S^2)$ 
the sequence $\{g^\psi_\eps\}\subset W^{1,1}((-1,1))$, defined by
\begin{gather*}
 g^\psi_\eps(s) := \int_{ \{u_\eps=s\} }\psi(y,
\niueps(y))~
\eps\vert\nabla u_\eps(y)\vert\,d\Ha^{2}(y),
\end{gather*}
converges strongly in $W^{1,1}((-1,1))$ to the function  $g^\psi(s):=\sqrt{2W(s)}V^0(\psi)$. 
Moreover, 
for $\LL^1$-almost every $s\in [-1,1]$ we have
\begin{equation}
\begin{split}
&\lim_{\eps\to 0^+}\var
\left(
\{u_\eps=s\},\star\niueps,\eps\vert\nabla 
u_\eps\vert
\right)=
\lim_{\eps\to 0^+}\var
\left(
\{u_\eps=s\},\star \niueps,\sqrt{2W(s)}
\right)
\\
=&\sqrt{2W(s)}V^0
\end{split}
\end{equation}
as oriented varifolds in $\Omega$.
\end{itemize}
\end{theorem}
\begin{remark}\label{alfffonso}
We can adapt the proof of Theorem \ref{theo:se2mni} to show
 that, under the weaker assumption that the hypothesis of Theorem \ref{the:RS} hold, the sequence $g_\eps^\psi$ converges strongly 
to $g^\psi$ in 
$W^{1,1}_{{\rm loc}}((-1,1))$ 
as $\eps\to 0^+$ for every $\psi\in C^1_c(\Omega)$.
\end{remark}
The next 
proposition shows  that a stronger convergence to zero of the 
discrepancies $\xiepsueps$ defined in \eqref{eq:def-xi-eps} holds, 
assuming the bounds in \eqref{eq:ener-bound}. Similar estimates have been obtained in \cite{PaTonz}, when $u_\eps$ is a local minimizer for $\mathpzc P_\eps$.
\begin{proposition}[Improved convergence of the 
discrepancies]\label{lem:conv-discr-BV}
Suppose that 
$\{u_\eps\} \subset \C^2(\Omega)$  is such that 
\eqref{eq:ener-bound} holds.
Then there exists a (not relabelled)
subsequence
such that
\begin{gather}
\nabla
\xi_{u_\eps}^\eps
~ \LL^3\rightharpoonup 
0\qquad \text{ as Radon measures on }\Omega,
\label{eq:discr-to-0-in-BV}
\\
\lim_{\eps\to 0^+}\|\xiepsueps\|_{L^p(\Omega)}=0 \qquad\text{ for every 
}p\in [1, 3/2).
\label{asburatto} 
\end{gather}
\end{proposition}
%
\section{Proof of Theorem \ref{the:Gamma-liminf-cond}}\label{sec:cavallona}

The present section is organized as follows. We start by proving two technical lemmata, 
namely Lemma \ref{lem:Kazumba}
and Lemma \ref{lem:conv-to-Hutch}. Then in  Section \ref{sec:staboia} 
we prove that  $V:=\lim_{\eps\to 0}\Vepsueps$ 
is a
 curvature varifold with generalized second fundamental
form in $L^2$, we show
 \eqref{eq:Aueps-to-AVmu-the} 
and  inequality \eqref{eq:Gamma-liminf-var}.

\begin{lemma}\label{lem:Kazumba}
Suppose that $\{u_\eps\} \subset \C^2(\Omega)$ is such that 
\begin{equation}\label{eq:granmatthias}
\sup_{0<\eps<1}\left\{\muepsueps(\Omega)+
\frac{1}{\eps}\int_\Omega \left\vert \eps\nabla^2u_\eps-
\frac{W^\prime(u_\eps)}{\eps}\niueps\otimes \niueps\right\vert^2\,dx
\right\}
<+\infty.
\end{equation}
Then there exists a (not relabelled) subsequence such that 
\begin{equation}\label{Kazumba}
\lim_{\eps\to 0^+}(V_{u_\eps}^\eps,R^\eps_{u_\eps})= (V,0)
\end{equation}
as measures function pairs on $G_2(\Omega)$ with values in $\R^3$,
where the varifold $V$ is defined in Theorem \ref{the:RS} (B). 
\end{lemma}
\proof
Since $\fepsueps=\mathrm{tr}\big(\eps\nabla^2u_\eps-\frac{1}{\eps}
W^\prime(u_\eps) \niueps\otimes\niueps \big)$, we have
\begin{gather*}
 \frac{1}{\eps}\int_\Omega (f_{u_\eps})^2\,dx\leq \frac{3}{\eps}\int_\Omega
 \left\vert \eps\nabla^2u_\eps-\frac{W^\prime(u_\eps)}{\eps}
\niueps\otimes \niueps\right\vert^2\,dx.
\end{gather*}
Hence, 
by \eqref{eq:granmatthias}, we can apply Theorem \ref{the:RS}, and select a 
(not relabelled) subsequence such that $\Vepsueps\to V$ as $\eps\to 0^+$ 
in the sense of varifolds, with
 $V=\var(\rectifiableset,\theta)\in\mathbf{IV}_2(\Omega)$.
Since on $\{\nabla u_\eps\neq 0\}$ we have
\begin{equation}\label{autobus}
R^\eps_{u_\eps} = \frac{\nabla\xiepsueps}{\eps\vert\nabla u_\eps\vert^2}=
B^\eps_{u_\eps}
\frac{\nabla u_\eps}{\vert\nabla u_\eps\vert},
\end{equation}
we conclude that
\begin{equation*}
\surftens \int\vert R^\eps_{u_\eps}\vert^2\,d\Vepsueps=
\int\left\vert\frac{\nabla\xiepsueps}{\eps\vert\nabla u_\eps\vert^2}
\right\vert^2\,d\tildemuepsueps\leq3
\int\vert B^\eps_{u_\eps}
\vert^2\,d\tildemuepsueps,
\end{equation*}
which is uniformly bounded with respect to $\eps$ in view of 
\eqref{eq:granmatthias}. 
By Theorem \ref{the:Hutch-comp-lsc} (i), 
we can select a further (not relabelled)
subsequence 
 such that $(\Vepsueps,R_{u_\eps}^\eps)$ 
converge  weakly as measure-function pairs on $G_2(\Omega)$ with 
values in $\R^3$ to $(V,R)$, for a certain $R\in L^2(V,\R^3)$.
In order to prove \eqref{Kazumba} we closely 
follow \cite[page 10]{Tonesff}. Let $\phi \in C^1_c(\Omega)$
and $R_i$ (resp. $R^\eps_{u_\eps,i}$) be the $i$-th
component of $R$ (resp. of $R^\eps_{u_\eps}$). 
 By \eqref{eq:the-RS-van-discr} we have
\begin{equation}\label{basso}
\surftens
\int R_i(x,S)\phi(x)\,dV(x,S)=\lim_{\eps\to 0^+}\int 
R^\eps_{u_\eps,i}\phi\,d\tildemuepsueps=
-\lim_{\eps\to 0^+}\int \partial_i\phi\,d\xiepsueps=0,
\end{equation}
where in the two last equalities we used
\eqref{autobus}, \eqref{grad:discr} and \eqref{eq:the-RS-van-discr} respectively. 

{}From \eqref{basso}, 
using that $\Vepsueps\to V=\var(\rectifiableset,\theta)\in\mathbf{IV}_2(\Omega)$ as varifolds, it follows
\begin{equation*}
\int R_i(x,S)\phi(x)\,dV(x,S)=0=
\int_M R_i(x,T_xM)\phi(x)\,\theta(x) d\Ha^2(x).
\end{equation*}
This implies that $R(x,T_xM)=0$ for 
$\mu_V=\theta\Ha^2\res M$-a.e. $x$, and \eqref{Kazumba} follows.
\endproof
\begin{remark}\rm
We need to consider $R^\eps_{u}$ as a function on $G_2(\Omega)$ and not just on $\Omega$ because $R^\eps_{u}$ appears in 
 the ``$\eps$-formulation''   of  \eqref{eq:def-Hutch-var} (see \eqref{primavera}), which   characterizes Hutchinson's curvature varifolds  via an ``integration by  parts'' formula involving test functions  in $C_c^1(G_2(\Omega))$.
\end{remark}
The following lemma shows that if \eqref{eq:granmatthias} holds then the 
varifold $V$ limit of the $V_{u_\eps}^\eps$ 
is a 
curvature varifold with generalized second fundamental
form in $L^2$.
\begin{lemma}\label{lem:conv-to-Hutch}
Suppose that \eqref{eq:granmatthias} holds. Then
\begin{equation}\label{eq:bound-ssfeps}
\sup_{0<\eps<1}\int\vert\sffeps\vert^2\,dV_{u_\eps}<\,
+\infty.
\end{equation}
Moreover the varifold $V$ in Lemma \ref{lem:Kazumba} 
is a curvature 
varifold with 
generalized second fundamental form ${\bf B}_V$ in $ L^2$
and, up to a subsequence,
\begin{align}
&\lim_{\eps\to 0^+}(V_{u_\eps},A^{u_\eps})=(V,A_{V}),
\label{eq:lim-Aeps}
\\
&\lim_{\eps\to 0^+}(V_{u_\eps},\sff_{u_\eps})=(V,\sff_{V}),
\label{eq:lim-Beps}
\end{align}
as measure-function pairs on $G_2(\Omega)$ with values in $\R^{3^3}$.
\end{lemma}

\proof
{}From the  definitions of $\sffeps$ 
and $B^\eps_{u_\eps}$ given in \eqref{eq:def-sff-lev-surf} 
and \eqref{eq:def-Aueps} respectively, 
we have
\begin{align}
\vert \sffeps\vert^2=&
\sum_{i,j,k=1}^3\left[\left(\frac{(\projueps)^T\nabla^2 u_\eps 
\projueps}
{\vert\nabla u_\eps\vert}\right)_{ij}\right]^2\left(\frac{\partial_k u_\eps}{\vert\nabla u_\eps\vert}\right)^2
\nonumber
\\
=&\sum_{i,j=1}^3\left[
\left(\frac{(\projueps)^T\nabla^2 u_\eps \projueps}{\vert\nabla 
u_\eps\vert}\right)_{ij}\right]^2=\left\vert\frac{(\projueps)^T\nabla^2 u_\eps \projueps}{\vert\nabla u_\eps\vert}\right\vert^2
\label{eq:acculo?}
\\ 
=&\left\vert\frac{(\projueps)^T\Big[\eps\nabla^2 u_\eps-
\frac{1}{\eps}W^\prime(u_\eps) \nabla 
u_\eps\otimes\nabla u_\eps/\vert\nabla u_\eps\vert^2\Big] 
\projueps}{\eps\vert\nabla u_\eps\vert}\right\vert^2
\leq\, \vert B^\eps_{u_\eps}\vert^2,
\nonumber
\end{align}
where in the last inequality we use \eqref{eq:Kill-Bill}. 
Integrating \eqref{eq:acculo?} with respect to $dV_{u_\eps}$ (see
\eqref{ehhh???} and \eqref{eq:def-xi-eps}) and using \eqref{eq:granmatthias}, 
we  conclude that \eqref{eq:bound-ssfeps} holds. 
Notice that by \eqref{eq:granmatthias}  the conclusions
of Theorem \ref{the:RS} hold.

By \eqref{def:Aijk} and \eqref{eq:granmatthias} we obtain also
\begin{equation*}
\sup_{0<\eps<1}\int \vert A^{u_\eps}\vert^2\,dV_{u_\eps}< +\infty.
\end{equation*}
This latter estimate together with $\sup_{0<\eps<1} 
\muepsueps(\Omega)< +\infty$, enables us to apply Theorem 
\ref{the:Hutch-comp-lsc} and 
conclude that, 
passing to a subsequence, there is $\widehat 
A \in L^2(V,R^{3^3})$ such that 
\begin{equation}\label{eq:conv-mfp-Aeps}
\lim_{\eps\to 0^+}(V_{u_\eps},A^{u_\eps})=(V,\widehat A)
\end{equation}
as measure-function pairs on $G_2(\Omega)$ with values on $R^{3^3}$.

Now we want to prove that actually $\widehat A(x,S)$ verifies equation \eqref{eq:def-Hutch-var} and hence that $V$ is a  curvature 
varifold with generalized second fundamental form in 
$L^2$, and 
 $\widehat A=A^{V}$. 
In doing this we closely follow \cite[Proposition 2]{Tonesff}.

Fix $1\leq i\leq 3$ and $\phi\in C_c^1(\Omega)$. 
Multiply equation \eqref{eq:def-chem-pot} by $\phi \partial_iu_{\eps}$. 
Integrating by parts we firstly obtain
\begin{equation}\label{eq:energytensor}
\int_\Omega\left[
\frac{\eps}{2}\vert\nabla u_\eps\vert^2 \partial_i\phi
-\eps \partial_iu_\eps \partial_j u_\eps \partial_j\phi
+\frac{W(u_\eps)}{\eps}\partial_i\phi\right]\, dx
=\int_\Omega \fepsueps ~ \phi~ \partial_iu_\eps\, dx.
\end{equation}
Hence
\begin{equation}\label{eq:trippa}
\int_\Omega
\left[\left(
\partial_i\phi-({\niueps})_{i} ({\niueps})_{j}\partial_j\phi
\right)\,\eps\vert\nabla u_\eps\vert^2
+
\phi\partial_i \xi^\eps_{u_\eps}
\right]~dx=\int_\Omega \fepsueps ~ \phi~ \partial_i 
u_\eps\,dx. 
\end{equation}
Let now $\varphi\in C^1_c(\Omega\times\R^{3\times 3})$, $\sigma>0$, and
define $\phi^\sigma \in C^1_c(\Omega)$ by
$$
\phi^\sigma(x)
:=\varphi\left(x,{\rm Id}-\frac{\nabla 
u_\eps(x)\otimes\nabla u_\eps(x)}{\sigma^2+\vert\nabla u_\eps(x)\vert^2}\right),
\qquad x \in \Omega.
$$
Using $\phi^\sigma$ in place of $\phi$ in \eqref{eq:trippa} and 
letting $\sigma\to 0^+$ we obtain
\begin{equation}\label{primavera}
\begin{split}
& \int_\Omega 
\left[
\projueps_{ij}\Big(\partial_j\varphi-\partial_j \left[({\niueps})_{l}
({\niueps})_k\right]
D_{m^{}_{lk}}\varphi\Big)-\frac{\fepsueps}{\eps\vert\nabla 
u_\eps\vert}\frac{\partial_i u_\eps}{\vert\nabla u_\eps\vert}
\varphi\right]\,d\widetilde\mu_{u_\eps}
\\
=&-
\int_\Omega\varphi\partial_i
\xi_{u_\eps}^\eps\,dx.
\end{split}
\end{equation} 
In \eqref{primavera} the integration is only on the subset of $\Omega$ 
where $\grad u_\eps \neq 0$, the function $\varphi$
is evaluated at $(x, {\rm Id} - \niueps(x) \otimes \niueps(x))$, and
$D_{m^{}_{lk}} \varphi$ is the derivative of $\varphi (x,\cdot)$ 
with respect to its $lk$-entry variable.
Next we notice that, by the definition of 
$\fepsueps$ and $A^{u_\eps}$ in \eqref{def:Aijk} we have
\begin{equation}\label{ferimaste}
\begin{split}
\frac{\fepsueps}{\eps\vert\nabla u_\eps\vert}\frac{\partial_i u_\eps}{\vert\nabla u_\eps\vert}=&
\mathrm{div}\left(\frac{\nabla u_\eps}{\vert\nabla u_\eps\vert}\right)\frac{\partial_i u_\eps}{\vert\nabla u_\eps\vert}
\\
& +\frac{1}{\eps\vert\nabla u_\eps\vert^2}\left[\frac{\eps\nabla^2u_\eps\nabla u_\eps\cdot\nabla u_\eps-\eps^{-1}W^\prime(u_\eps)\vert\nabla u_\eps\vert^2}{\vert\nabla u_\eps\vert^2}\right]\partial_i u_\eps
\\
=&
A^{u_\eps}_{jij}(x, \projueps
)+\frac{1}{\eps\vert\nabla u_\eps\vert^2}
\left(\nabla \xiepsueps \niueps\otimes\niueps\right)_i.
\end{split}
\end{equation}
Inserting \eqref{ferimaste} into \eqref{primavera},
and recalling the definition of $\Vepsueps$, $A^{u_\eps}$ and 
$R^\eps_{u_\eps}$ given in \eqref{ehhh???}, \eqref{def:Aijk}
and \eqref{erbarmedich} respectively, we have that 
equality \eqref{primavera} becomes 
\begin{equation*}
\begin{split}
&\int(S_{ij}\partial_j\varphi+A^{u_\eps}_{ijk}D_{m_{jk}}
\varphi-A_{jij}^{u_\eps}\varphi)\,d\Vepsueps(x,S)
\\
=&-
\int (R^\eps_{u_\eps}(x,S) ~S)_i \varphi(x,S)\,d\Vepsueps(x,S),
\end{split}
\end{equation*}
where $\varphi$ on the left hand side is evaluated at $(x,S)$.
Passing to the limit as $\eps\to 0^+$, by the convergence of 
$\{\Vepsueps\}$ to $V$,
 \eqref{eq:conv-mfp-Aeps} and Lemma \ref{lem:Kazumba}, we get
\begin{equation*}
\int\left
(S_{ij}\partial_j\varphi+\widehat{A}_{ijk}D_{m_{jk}}\varphi-\widehat A_{jij}\varphi\right)\,dV(x,S)=0,
\end{equation*}
that is $V$ is a curvature varifold with  generalized second
fundamental form in $L^2$, and $A^{V}=\widehat A$.

In order to get 
\eqref{eq:lim-Beps} we proceed as follows. 
Let $V=\var(\rectifiableset
,\theta)$. We define 
\begin{align*}
&\overline{\projueps}:G_2(\Omega)\to \R^{3\times 3},\qquad (x,S)\to 
\projueps(x),
\\
&\overline{P^V}:G_2(\Omega)\to \R^{3\times 3},\qquad (x,S)\to P^{\mathcal M}(x),
\end{align*}
where $P^{\rectifiableset
}(x)$ is the orthogonal 
projection matrix of $\R^3$ 
onto the tangent plane $T_x \rectifiableset
\in G_{2,3}$ to $\rectifiableset$ at $x$ 
(recall that $T_x \rectifiableset
$ is well defined $\Ha^2\res \rectifiableset
$-almost
everywhere by the $2$-rectifiability of $\rectifiableset
$, see \cite{AFPa00}).
By Remark \ref{rem:lerchiopattino}
we have that  
the convergence of $\Vepsueps$ to $V$ as varifolds implies that 
$(\Vepsueps,\overline{\projueps})\to (V,\overline{P^V})$ as $\eps\to 0^+$ in the  $L^2$-strong convergence as measure-function pairs on $G_2(\Omega)$ with values in $\R^{3 
\times 3
}$. Hence, by  \eqref{Bjork} and Lemma \ref{lem:Moser} 
we obtain \eqref{eq:lim-Beps}.
\endproof

Note that the left hand side of \eqref{eq:energytensor} can also be written as
$\int_\Omega T_\eps^{ij} \partial_j \phi~dx$, where $T_\eps^{ij}$
is the so-called energy-momentum tensor, defined as 
$T_\eps^{ij} := \left(\frac{\eps}{2} \vert \grad u\vert^2 
+ \frac{1}{\eps} W(u)\right)\delta_{ij} - \eps \partial_i u \, \partial_j u$.

\subsection{Proof of \eqref{eq:Gamma-liminf-var}}\label{sec:staboia}
{}From the definition of $\mathpzc W_\eps$  
in \eqref{veroe'benpindemonte} we have
\begin{equation}\label{pallissime}
\begin{split}
\mathpzc W_\eps(u_\eps)=
-\frac{\elGauss}{2\eps}\int_\Omega\left\vert\eps\nabla^2 
u_\eps-\frac{W^\prime(u_\eps)}{\eps}\niueps\otimes\niueps \right\vert^2\,dx
+\frac{\elbend+\elGauss}{2\eps}\int_\Omega (\fepsueps)^2~dx.
\end{split}
\end{equation}
{}From \eqref{eq:constr}, 
\eqref{eq:ener-bound} and \eqref{pallissime}
it follows
that \eqref{eq:granmatthias} holds. Hence by Lemma \ref{lem:conv-to-Hutch}
 we can conclude that $V$ is a curvature varifold with generalized 
second fundamental form ${\bf B}_V$ in $L^2$, and $A_{V}\in L^2(\mu_V)$ 
and also that \eqref{eq:Aueps-to-AVmu-the} is verified.
In order to prove the $\Gamma-\liminf$ inequality \eqref{eq:Gamma-liminf-var} 
we observe that, by \eqref{eq:acculo?}, we have
\begin{equation}\label{eq:abbestia?}
\begin{split}
\mathpzc W_\eps(u_\eps)\geq&
\int \left[\frac{-\elGauss}{2}
\vert\sffeps
\vert^2+
\frac{\elbend+\elGauss}{2}(\Hepsueps)^2
\right]\,d\tildemuepsueps
\\
=&\surftens\int \frac{-\elGauss}{2}
\vert\sffeps(x, \projueps)
\vert^2\,d\Vepsueps+\int
\frac{\elbend+\elGauss}{2} (\Hepsueps)^2
\,d\tildemuepsueps.
\end{split}
\end{equation}
By \eqref{eq:abbestia?}, \eqref{eq:lim-Beps}, 
and Theorem \ref{the:Hutch-comp-lsc}, we have
\begin{equation*}
\begin{split}
\liminf_{\eps\to 0^+}\,\mathpzc W_\eps(u_\eps)\geq&\surftens
\liminf_{\eps\to 0^+}\int \frac{-\elGauss}{2}\vert
\sff_{u_\eps}\vert^2\,d\Vepsueps+\liminf_{\eps\to 0^+}\int
\frac{\elbend+\elGauss}{2}(\Hepsueps)^2
\,d\tildemuepsueps
\\
\geq &
\surftens \int\left[ \frac{\elbend}{2}
\vert\Hv_{V}\vert^2+\frac{\elGauss}{2}\left(
\vert \Hv_{V}\vert^2-\vert
\sff_{V}\vert^2\right)\right]\,dV,
\end{split}
\end{equation*}
which proves \eqref{eq:Gamma-liminf-var}.

\section{Proofs of Theorem \ref{the:Gamma-limsup}
and of 
Corollary \ref{cor:Gamma-conv-smooth-points}}
\label{prendodimiraquestalberoinmezzoalgiardino}%
We prove Theorem \ref{the:Gamma-limsup} in the case $\Omega=
\R^3$. The case of a bounded $\Omega$ can be proved almost in the same way.

We will construct a sequence $\{u_\eps\}\subset H^2(\R^3)$ satisfying the 
thesis. 
To conclude the proof it is enough to mollify each $u_\eps$ and use a 
standard diagonal argument to obtain a new sequence $\{\widehat 
u_\eps\}\subset 
C^2(\R^3)$ still satisfying 
\eqref{eq:L1-conv-Gammalimsup}, \eqref{eq:Rm-conv-Gammalimsup}, \eqref{eq:Gammalimsup}.
 
 We consider $u_\eps\in H^2(\R^3)$ as in \cite{BePa}. 
Let $d(\cdot)$ be the signed distance function from
$\partial E$, as defined in Proposition \ref{prop:dist-VS-curv}, and let $\gamma(s):=\tanh(s)$.
For any $0 < \eps < 1$ and $s \in \R$, let $\gamma_\eps(s) :=
\gamma(s/\eps)$ and $\tgae$ be defined as follows:
$\tgae
 := \gae$ in $(0,\eps|\log \eps|)$,
$\tgae :=
p_\eps$ in $(\eps|\log \eps|, s^0_\eps)$,
$\tgae :=
+1$ in $(s^0_\eps, +\infty)$, and
$\tgae(s) :=
-\tgae(-s)$ if $s < 0$.
Here, $p_\eps$ is an arc of parabola on $(\eps \vert
\log\eps\vert, s^0_\eps)$ connecting the points $(\eps|\log \eps|,
\gae(\eps |\log \eps|))$ and $(s^0_\eps,1)$, that is $p_\eps(s) :=
-a_\eps(s-s^0_\eps)^2 + 1$, $a_\eps >0$. To find $a_\eps$ and
$s^0_\eps$, we impose the condition $\tgae \in H^2(\R)$,
that gives
$s^0_\eps = \eps + \eps^3 + \eps |\log \eps|$ and $a_\eps =
\frac{2}{(1 + \eps^2)^3}$.

We define
\begin{equation}\label{ve}
u_\eps(x) := \tgae(d(x)).
\end{equation}
Then \eqref{eq:L1-conv-Gammalimsup} and 
\eqref{eq:Rm-conv-Gammalimsup} follow directly from \cite{BePa}, and it remains to prove only \eqref{eq:Gammalimsup}.

To this aim we notice that, since
$\nabla^2 u_\eps=\tgae^\prime(d)\nabla^2 d+\tgae^{\prime\prime}(d)
\nabla d\otimes\nabla d$, we have
\begin{itemize}
\item[-]
 in $U_\eps:=\{-\eps\vert\log\eps\vert<d(x)<\eps\vert\log\eps\vert\}$
\begin{align}
&
B^\eps_{u_\eps}
=\frac{\gamma^\prime(d/\eps)\nabla^2 d+\eps^{-1}\Big(\gamma^{\prime\prime}
(d/\eps)
-W^\prime(\gamma(d/\eps))\Big)\nabla d\otimes\nabla d}{\vert\gamma^\prime(d/\eps)\vert}
=\nabla^2 d,
\label{eq:Aueps-Ueps}
\\
&\notag
\\
&\Hepsueps=\Delta d;
\label{eq:Heps-Ueps}
\\
&\notag
\end{align}
\item[-] in $\mathcal V_\eps
:=\{\eps\vert\log\eps\vert<\vert d(x)\vert<  s^0_\eps\}$
\begin{align}
&
B^\eps_{u_\eps}
=\nabla^2 d+\frac{1}{\eps p_\eps^\prime(d)}\left(\eps p_\eps^{\prime\prime}(d)-
\frac{W^\prime(p_\eps(d))}{\eps}\right)\nabla d\otimes\nabla d,
\label{eq:Aueps-Veps}
\\
&\Hepsueps=\Delta d +\frac{1}{\eps p_\eps^\prime(d)}\left(\eps p_\eps^{\prime\prime}(d)-
\frac{W^\prime(p_\eps(d))}{\eps}\right).
\label{eq:Hueps-Veps}
\end{align}
\end{itemize}
Let us now derive some estimates in $\mathcal V_\eps$.
Let $x \in \mathcal V_\eps$; 
then $1 \geq \displaystyle u_\eps(x)
\geq p_\eps(\eps |\log \eps|) = 1 - \frac{2 \eps^2}{1 + \eps^2}$.
Hence $\displaystyle |W^\prime(u_\eps(x))| =
|4u_\eps(x)(1-u_\eps(x))(1+u_\eps(x))| \leq \frac{16\eps^2}{1 + \eps^2}$,
so that $\eps^{-1} W^\prime(u_\eps) = O(\eps)$. Moreover
$\eps p_\eps^{\prime\prime}(d) = O(\eps)$,
so that
\begin{equation}\label{eq:pippo}
-\eps p_\eps''(d) +
 \frac{W^\prime(u_\eps)}{\eps} = O(\eps).
\end{equation}
Moreover since  $\eps|p_\eps'(s)|^2 = \frac{ 8\eps(s - \eps - \eps^3 -
\eps|\log \eps|)^2 }{ (1+\eps^2)^6 }$, making the change of
variable $\sigma = s- \eps|\log \eps|$, it follows
\begin{equation}\label{eq:pazzo}
\int_{\eps |\log \eps|}^{\eps + \eps^3 + \eps |\log \eps|}
\eps|p_\eps'(s)|^2~ds = \frac{32 \eps}{(1 + \eps^2)^6}
\int_0^{\eps + \eps^3} (\tau - \eps - \eps^3)^2~d\tau =
O(\eps^4),
\end{equation}
as $\eps \to 0^+$

By \cite{BePa} it follows that
\begin{equation}\label{eq:Heps-Gammalimsup}
\lim_{\eps\to 0^+}\mathpzc H_\eps(u_\eps)=\surftens\int_{\partial E}(H_{\partial 
E})^2\,d\Ha^2.
\end{equation}
Eventually  we have 
\begin{equation}\label{eq:Keps-Gammalimsup} 
\begin{split}
\lim_{\eps\to 0^+}\mathpzc{K}_\eps(u_\eps)=&
\lim_{\eps\to 0^+} 
\Bigg\{\int_{U_\eps}
\sum_{1\leq i<j\leq 
j}
\det([B^\eps_{u_\eps}
]_{ij})
\,\eps\vert\nabla u_\eps\vert^2dx 
+\int_{\mathcal V_\eps}\frac{(\Hepsueps)^2-\vert B^\eps_{u_\eps}
\vert^2}{2}\,\eps\vert\nabla u_\eps\vert^2dx \Bigg\}
\\
=&\lim_{\eps\to 0^+}
\Bigg\{
\int_{U_\eps}\sum_{1\leq i<j\leq 3}\det([\nabla^2 
d]_{ij})\frac{1}{\eps}\vert\gamma^\prime(d/\eps)\vert^2\, dx 
\\
&+\frac{1}{2\eps}\int_{\mathcal V_\eps}
\left[\eps p_\eps^\prime\Delta d +
\left(\eps p_\eps^{\prime\prime}(d)-\frac{W^\prime(p_\eps(d))}{\eps}\right)\right]^2\,dx
\\
&-\frac{1}{2\eps}\int_{\mathcal V_\eps}
\left\vert\eps p_\eps^\prime\nabla^2 d +
\left(\eps p_\eps^{\prime\prime}(d)-
\frac{W^\prime(p_\eps(d))}{\eps}\right)\nabla d\otimes\nabla d\right\vert^2\,dx
\Bigg\}
\\
=&\lim_{\eps\to 0^+}\left(
\int_{U_\eps}\sum_{1\leq i<j\leq 3}\det([\nabla^2 
d]_{ij})\frac{1}{\eps}\vert\gamma^\prime(d/\eps)\vert^2\, dx+O(\eps)\right)
\\
=&\surftens\int_{\partial E} K_{\partial E}\,d\Ha^2,
\end{split}
\end{equation}
where in the last equality we use Proposition \ref{prop:dist-VS-curv}. 
Hence, by \eqref{eq:Heps-Gammalimsup} and 
\eqref{eq:Keps-Gammalimsup} we deduce that  \eqref{eq:Gammalimsup} holds.

\subsection{Proof 
of 
Corollary \ref{cor:Gamma-conv-smooth-points}}\label{ziocantante}
If  
$E$ has smooth boundary in $\Omega$, as in the proof of Theorem \ref{the:Helf-lsc},  we can use the locality
 of the 
generalized second fundamental form for Hutchinson's curvature varifolds (see \cite{Mant}) together with
\begin{equation*}
\surftens \Ha^2 \res\partial E\leq  \mu=\surftens\mu_V
\qquad \text{as Radon~measures,}
\end{equation*}
to conclude that
\begin{equation*}
\surftens \int\left[ \frac{\elbend}{2}\vert\Hv_{V}\vert^2+
\frac{\elGauss}{2}\left(
\vert\Hv_{V}\vert^2-\vert\sff_{V}\vert^2\right)\right]\,dV
\geq \surftens\Helf(E).
\end{equation*}
The thesis is then a direct consequence of Theorems \ref{the:Gamma-liminf-cond}, \ref{the:Gamma-limsup} and \ref{the:Helf-lsc}.
\section{Proof of 
Theorem \ref{theo:se2mni}}
\label{maanvediunpodanna}

Firstly we notice that we can assume (up to  selecting  a subsequence) that  $\Vepsueps$ 
converge as varifolds to the curvature 
varifold $V\in\mathbf{IV}_2(\Omega)$ and that \eqref{eq:Aueps-to-AVmu-the} holds. Moreover, since
$\Vzeroepsueps(G^0_2(\Omega))=\muepsueps(\Omega)$, by 
\eqref{eq:ener-bound}, we can 
extract a further subsequence such that 
$\Vzeroepsueps$ converge as Radon measures  to a Radon measure $V^0$ on 
$G^0_2(\Omega)$, and also that $\q_\sharp  V^0=V$
(notice that for the moment $V^0$ is rectifiable but not necessarily integral). 
Eventually, without loss of generality, we can also assume that 
 $$
 \liminf_{\eps\to 0^+} 
\mathpzc W_\eps(u_\eps)=\lim_{\eps\to 0^+}\mathpzc W_\eps(u_\eps)<+\infty.
 $$
The present section is organized as follows.
We firstly prove Lemma \ref{tuttoperCicciobelloRicchione}, from which
Theorem \ref{theo:se2mni}-(B)  follows. Then, in Proposition \ref{smandruppone}, we conclude the proof of Theorem \ref{theo:se2mni}-(A) showing that $V^0\in\mathbf{IV}^0_2(\Omega)$.
\begin{lemma}\label{tuttoperCicciobelloRicchione}
Let $u_\eps\in C^2(\Omega)$ be such that \eqref{eq:granmatthias} holds
and $\liminf_{\eps \to 0^+} \muepsueps(\Omega)>0$. Suppose $V^0$ is such that
$\lim_{\eps\to 0^+}\Vzeroepsueps
=V^0$ as oriented varifolds. Then 
there exists a
(not relabelled) subsequence of $\{u_\eps\}$ such that for $\LL^1$-almost every $s\in [-1,1]$ we have
\begin{equation}
\begin{split}
&\lim_{\eps\to 0^+}\var
\left(
\{u_\eps=s\},\star\niueps,\eps\vert\nabla 
u_\eps\vert
\right)=
\lim_{\eps\to 0^+}\var
\left(
\{u_\eps=s\},\star \niueps,\sqrt{2W(s)}
\right)
\\
=&\sqrt{2W(s)}V^0
\end{split}
\end{equation}
as oriented varifolds on $\Omega$. 
\end{lemma}
\begin{remark}\label{rem:heidi-pornostar}
When $\liminf_{\eps\to 0^+}\muepsueps(\Omega)>0$, by \cite[Lemma 4.4]{Chen} (see also \cite[Proposition 3.4]{RS}) we can conclude that, 
up to a subsequence, $\{u_\eps=s\}\neq\emptyset$ for every $s\in (-1,1)$.
\end{remark}

\proof
Let us firstly remark that 
on one hand for $\psi\in C_c^0(\Omega\times 
\mathbb S^2)$  we can define 
$\psi^*\in C_c^0(G_2^0(\Omega))$ as
$\psi^\star(x,\tau):=\psi(x,\nu^\tau)$.  
On the other hand for  
$\phi\in C^0_c(G^0_2(\Omega))$
 we can define $\phi_\star\in  
C_c^0(\Omega\times \mathbb S^2)$ as $\phi_\star(x,\xi)
:=\phi(x,\star\xi)$. This means that the convergence as oriented varifolds of $\var(\{u_\eps=s\},\star\niueps,1)$ is equivalent to the convergence of $\Ha^2\res\{u_\eps=s\}\otimes\delta_{\niueps}$ as 
measures on $\Omega\times\mathbb S^2$.
Moreover for a 
given $\psi\in
 C_c^1(\Omega\times \mathbb S^2)$ we can find $\overline\psi\in C^1_c(\Omega\times\R^3)$ 
such that $\psi(x,\xi)=\overline\psi(x,\xi)$ for every 
$\xi\in \mathbb S^2$, and $\|\overline\psi\|_{L^{\infty}(\Omega\times\R^3)}\leq\|\psi\|_{L^{\infty}(\Omega\times\mathbb S^2)}$.

Let $\psi\in C^1_c(\Omega\times \mathbb S^2)$ and define
$g^\psi_{\eps}:\R\,\to\,[0,+\infty)$ as in the statement, i.e.,
$$
g^\psi_\eps(s) := \int_{ \{u_\eps=s\} }\psi(y,
\niueps(y))~
\eps\vert\nabla u_\eps(y)\vert\,d\Ha^{2}(y).
$$
We extend $\psi$ to a function of class $C^1_c(\Omega\times B)$, where
$B:= \{\xi \in \R^3 : \frac{1}{2} < \vert \xi\vert < 2\}$, 
and we still denote by $\psi = \psi(x,\xi)$ such an extension.
Fixed $\delta\in (0,1/2]$ we set $I_\delta:=[-1+\delta,1-\delta]$.
Let $\eta\in C^\infty_c(I_\delta)$.
For fixed $\eps>0$ and $\sigma \neq 0$, 
we define $\psi^\sigma\in C^1_c(\Omega\times\R^3)$  as
\begin{equation*}
\psi^\sigma(x):=\psi\left(x,\frac{\nabla u_\eps(x)}{\sigma^2+\vert\nabla u_\eps(x)\vert}\right),
\end{equation*}
so that, since
 $\psi\in C^1_c(\Omega\times B)$, we obtain $\psi^\sigma\equiv 0$ on $\{\nabla u_\eps=0\}$.
We then have, 
using the coarea formula, 
\begin{align*}
&\int_\R\eta^\prime g^{\psi^\sigma}_{\eps}\,ds=
\int_\Omega \eps \eta^\prime(u_\eps)\psi^\sigma\vert\nabla u_\eps\vert^2\,dx
=\int_\Omega \eps \psi^\sigma\nabla (\eta(u_\eps))\cdot\nabla u_\eps\,dx
\\
=&-\int_\Omega \eps
\psi^\sigma\eta(u_\eps) \Delta u_\eps\,dx
-\int_\Omega\eps\eta(u_\eps)\nabla \psi^\sigma\cdot\nabla u_\eps\,dx.
\end{align*}

Letting $\sigma\to 0$ we obtain
\begin{equation}\label{eq:ponzo}
\begin{split}
&\int_\R\eta^\prime g^\psi_{\eps} \,ds
=-\int_{\Omega_\eps} \eps \eta(u_\eps)\
\psi \Delta u_\eps\,dx
\\
-&\int_{\Omega_\eps}\eps\eta(u_\eps)\nabla\psi\cdot\nabla u_\eps\,dx
-\int_{\Omega_\eps}\eps\eta(u_\eps)D_{\xi_j}\psi(x,
\niueps
)\partial_k (\nu^{u_\eps})_j\partial_ku_\eps\,dx,
\end{split}
\end{equation}
where $\Omega_\eps:=\Omega\cap\{\nabla u_\eps\neq 0\}$.

Adding and subtracting the term
$\int_{\Omega_\eps}\eta(u_\eps)\psi\frac{W^\prime(u_\eps)}{\eps}\,dx$,
 observing that the last addendum on the right hand side
of \eqref{eq:ponzo} can be written as 
$$
-\int_{\Omega_\eps}\eta(u_\eps)D_\xi\psi(x,\niueps)\projueps
\eps\nabla^2u_\eps
\frac{\nabla u_\eps}{\vert\nabla u_\eps\vert}\,dx,
$$
and since $P^{u_\eps} \nu_{u_\eps}\otimes \nu_{u_\eps}=0$, from 
\eqref{eq:ponzo} we obtain
\begin{align}
&\int_\R\eta^\prime g^\psi_{\eps} \,ds=
\int_{\Omega_\eps}\eta(u_\eps)\psi\left(-\eps\Delta u_\eps+\frac{W^\prime(u_\eps)}{\eps}\right)\,dx
-\int_{\Omega_\eps}\eps\eta(u_\eps)\nabla\psi\cdot\nabla u_\eps\,dx
\notag
\\
-&\int_{\Omega_\eps}\eta(u_\eps)D_\xi\psi(x,\niueps)\Big(\projueps
\Big(\eps\nabla^2u_\eps-
\frac{W^\prime(u_\eps)}{\eps}
\niueps
\otimes\niueps
\Big)
\Big)
\frac{\nabla u_\eps}{\vert\nabla u_\eps\vert}\,dx
\label{eq:ponzoreprise}
\\
-&\int_{\Omega_\eps}\eta(u_\eps)\psi\frac{W^\prime(u_\eps)}{\eps}\,dx.
\notag
\end{align}

Since for every $t\in I_\delta$ we have
%
$\vert W^\prime(t)\vert=\vert t(1-t^2)\vert
\leq \frac{4(1-\delta)}{\delta} W(t)$,
%
we can conclude that
\begin{equation*}
\begin{aligned}
\left\vert\int_\R \eta^\prime g^\psi_{\eps} \,ds\right\vert
\leq & \|\eta\|_{L^\infty(I_\delta)}
\|\psi\|_{L^\infty(\Omega\times 
\mathbb S^2)}\|\fepsueps\|_{L^1(\Omega)}
\\
& +\eps^{1/2}\|\eta\|_{L^\infty(I_\delta)}\|\nabla\psi\|_{L^\infty(\Omega\times \mathbb S^2)}\left(\int_\Omega\eps\vert\nabla u_\eps\vert^2\,dx\right)^{1/2}
\\
& +\|\eta\|_{L^\infty(I_\delta)}\|D_\xi\psi\|_{L^\infty(\Omega\times 
\mathbb S^2)}
\int_\Omega\Big\vert\eps\nabla^2 u_\eps-\frac{W^\prime(u_\eps)}{\eps}
\niueps
\otimes
\niueps
\Big\vert\,dx
\\
& + \|\eta\|_{L^\infty(I_\delta)}\|\psi\|_{L^\infty(\Omega\times \mathbb S^2)}\frac{4(1-\delta)}{\delta}\int_\Omega\frac{W(u_\eps)}{\eps}\,dx.
\end{aligned}
\end{equation*}

{}From this inequality we can deduce
that there exists $g^\psi\in BV_{{\rm loc}}([-1,1])$ such that 
$g^\psi_\eps
\to g^\psi$ in $L^1_{{\rm loc}}([-1,1])$ and $\LL^1$-almost everywhere in $[-1,1]$.

Next, for any fixed $\psi\in C^1_c(\Omega)$, we consider 
the functions $\widehat g_\eps^\psi:\R\,\to\,[0,+\infty)$ defined as
$$
\widehat g_\eps^\psi(s) := 
\sqrt{2W(s)}\int_{\{u_\eps=s\}}\psi(y,
\niueps
(y))\,d\Ha^{2}(y),
$$
and we claim that as $\eps\to 0^+$ the sequence $\{\widehat g_\eps^\psi\}$ 
converges in $L^1_{{\rm loc}}([-1,1])$ and $\mathcal L^1$-almost everywhere 
to $g^\psi$.
In order to prove the claim, let $\delta>0$. By \eqref{eq:lim-tuttuguale}
we have
\begin{align*}
&\lim_{\eps\to 0^+}\int_{I_\delta}\Big\vert\widehat 
g_\eps^\psi -g^\psi \Big\vert\,ds\leq\lim_{\eps\to 
0^+}\left(
\int_{I_\delta}\Big\vert\widehat 
g^\psi_\eps -g^\psi_\eps \Big\vert\,ds+
\int_{I_\delta}\vert g_\eps^\psi -g^\psi \vert\,ds \right)
\\
= &\lim_{\eps\to 0^+}
\left(\int_{I_\delta}\Big\vert\int_{\{u_\eps=s\}} 
\psi(\sqrt{2W(s)}-\eps\vert\nabla u_\eps\vert)\,d\Ha^2\Big\vert\,ds+
O(\eps)\right)
\\
\leq&\lim_{\eps\to 0^+}\int_{I_\delta}\int_{\{u_\eps=s\}} 
\Big\vert\psi(\sqrt{2W(s)}-\eps\vert\nabla u_\eps\vert)\Big\vert\,d\Ha^2\,ds
\\
\leq &  2\|\psi\|_{L^\infty(\Omega\times \mathbb S^2)}\lim_{\eps\to 
0^+}\int_{\Omega\cap\{u_\eps\in I_\delta\}}\left\vert
\sqrt{\frac{W(u_\eps)}{\eps}}-\sqrt{\frac{\eps}{2}}\vert\nabla u_\eps\vert\right\vert\,\sqrt{\frac{\eps}{2}}\vert\nabla u_\eps\vert\,dx
\\
\leq& 2\|\psi\|_{L^\infty(\Omega\times \mathbb S^2)}\lim_{\eps\to 
0^+}\int_{\Omega}\left\vert
\sqrt{\frac{W(u_\eps)}{\eps}}-\sqrt{\frac{\eps}{2}}\vert\nabla u_\eps\vert\right\vert\,\left(\sqrt{\frac{\eps}{2}}\vert\nabla u_\eps\vert\,+\sqrt{\frac{W(u_\eps)}{\eps}}\right)\,dx
\\
=& 2\|\psi\|_{L^\infty(\Omega\times \mathbb S^2)}\lim_{\eps\to 
0^+}\int_{\Omega}\vert
\xiepsueps\vert\,dx
=0,
\end{align*}
which shows the claim.
Since on $I_\delta$ we have $\sqrt{(2\delta-\delta^2)/2}\leq\sqrt{2W(s)}\leq\sqrt{2}$,
we can also conclude that the 
sequence of functions
 \begin{gather*}
 h^\psi_\eps
:\R\,\to\,[0,+\infty),
 \qquad
 h^\psi_\eps(s) :=
\frac{\widehat g^\psi_\eps(s)}{\sqrt{2W(s)}}=\int_{\{u_\eps=s\}}\psi(y,
\niueps
(y))d\Ha^{2}(y),
 \end{gather*}
 is equibounded in $L^1_{{\rm loc}}([-1,1])$ 
and converges in $L^1_{{\rm loc}}([-1,1])$ to 
 %
 %
\begin{equation}\label{h_psi}
h^\psi =\frac{ g^\psi}{\sqrt{2W}}.
\end{equation}
%
%
Next we refine formula \eqref{eq:ponzo}, by proving that, for every $\delta>0$, every $\psi\in C^1_c(\Omega)$ and $\eta\in C^\infty_c(I_\delta)$, we have 
\begin{equation}\label{eq:constipation}
\lim_{\eps\to 0^+}\int_{I_\delta}\eta^\prime g^\psi_\eps \,ds=
\int_{I_\delta}\eta \Big(\frac{d}{ds}\sqrt{2W}\Big) h^\psi\,ds.
\end{equation}
To this aim we start noticing that
\begin{align*}
&\left\vert\int_{\Omega_\eps}\eta(u_\eps)\psi\frac{W^\prime(u_\eps)}{\eps}\,dx- \int_{\Omega}\eta(u_\eps)\psi\frac{W^\prime(u_\eps)}{\sqrt{2W(u_\eps)}}\vert\nabla u_\eps\vert\,dx\right\vert
\\
\leq&\|\eta\|_{L^\infty(I_\delta)}\|\psi\|_{L^\infty(\Omega\times 
\mathbb S^2)}
\int_{\Omega\cap\{u_\eps\in I_\delta\}}\frac{\vert W^\prime(u_\eps)\vert}{\eps^{1/2}\sqrt{W(u_\eps)}}\left\vert\sqrt{\frac{W(u_\eps)}{\eps}}-\sqrt{\frac{\eps}{2}}\vert\nabla u_\eps\vert\right\vert\,dx
\\
\leq&\|\eta\|_{L^\infty(I_\delta)}\|\psi\|_{L^\infty(\Omega\times \mathbb S^2)}\frac{4(1-\delta)}{\delta}
\left(\int_{\Omega}
\frac{W(u_\eps)}{\eps}\,dx\right)^{1/2}
\left\|\sqrt{\frac{W(u_\eps)}{\eps}}-\sqrt{\frac{\eps}{2}}\vert\nabla u_\eps\vert\right\|_{L^2(\Omega)},
\end{align*}
which, by \eqref{eq:the-RS-van-discr}, vanishes as $\eps\to 0^+$. Then, by the $L^1(I_\delta)$ convergence
 of $h^\psi_\eps$,
the coarea formula and the Lebesgue's 
Dominated Convergence theorem, we have
\begin{align*}
&\lim_{\eps\to 0^+}\int_{\Omega_\eps}\eta(u_\eps)\psi\frac{W^\prime(u_\eps)}{\eps}\,dx=
\lim_{\eps\to 0^+}\int_{\Omega\cap\{u_\eps\in I_\delta\}}
\eta(u_\eps)\psi\frac{W^\prime(u_\eps)}{\sqrt{2W(u_\eps)}}\vert\nabla 
u_\eps\vert\,dx
\\
=&\lim_{\eps\to 0^+}
\int_{I_\delta}\eta \left(\frac{d}{ds}\sqrt{2W}\right)h^\psi_\eps \,ds
=\int_{I_\delta}\lim_{\eps\to 
0^+}\left(\eta\left(\frac{d}{ds}\sqrt{2W}\right)h^\psi_\eps\right)\,ds
\\
=&\int_{I_\delta}\eta\left(\frac{d}{ds}\sqrt{2W}\right)h^\psi\,ds.
\end{align*}
In order to obtain \eqref{eq:constipation} it is then enough to 
plug the  following estimates in \eqref{eq:ponzoreprise}:
\begin{gather*}
\left\vert\int_\Omega\eta(u_\eps)\psi \fepsueps\,dx\right\vert
\leq\eps^{1/2}\|\eta\|_{L^\infty(I_\delta)}\|\psi\|_{L^\infty(\Omega\times \mathbb S^2)}
\sqrt{\LL^n(\Omega)}\|\eps^{-1} \fepsueps\|_{L^2(\Omega)},
\end{gather*}
\begin{gather*}
\left\vert\int_\Omega\eps\eta(u_\eps)\nabla\psi\cdot\nabla u_\eps\,dx
\right\vert\leq\eps^{1/2} \|\eta\|_{L^\infty(I_\delta)}\|\nabla\psi
\|_{L^{\infty}(\Omega\times \mathbb S^2)}
\sqrt{\LL^n(\Omega)}\left(\int_{\Omega}
\eps\vert\nabla u_\eps\vert^2\,dx\right)^{1/2},
\end{gather*}
and
\begin{align*}
&\left\vert\int_\Omega\eta(u_\eps)D_\xi\psi(x,\niueps)\left(
\projueps
\left(\eps\nabla^2u_\eps-\frac{W^\prime(u_\eps)}{\eps}
\niueps
\otimes\niueps
\right)\right)\frac{\nabla u_\eps}{\vert\nabla u_\eps\vert}\,dx\right\vert
\\
\leq&\|\eta\|_{L^\infty(I_\delta)}\|D_\xi\psi\|_{L^\infty(\Omega\times 
\mathbb S^2)}\eps^{1/2}\|B^\eps_{u_\eps}\|_{L^2(\tildemuepsueps)}.
\end{align*}

We are now in a
position to prove that the distributional derivative 
of the function $h^\psi$ in \eqref{h_psi} 
is zero in $I_\delta$. 
In fact by \eqref{eq:constipation}, the definition of $h^\psi_\eps$ and Lebesgue's 
Dominated Convergence Theorem we have
\begin{align*}
&\int_{I_\delta}\eta' \sqrt{2W}h^\psi \,ds=
\lim_{\eps\to 0^+}\int_{I_\delta}g^\psi_\eps \,ds
=-\int_{I_\delta}\eta
\left(\frac{d}{ds}\sqrt{2W}\right)h^\psi \,ds,
\end{align*}
that is, for every $\eta\in C^\infty_c(I_\delta)$ we have
\begin{gather}\label{eq:sjh}
\int_{I_\delta}\frac{d}{ds}\left(\eta \sqrt{2W}\right)h^\psi \,ds=0.
\end{gather}
Since $\sqrt{2W}\geq \sqrt{\frac{2\delta-\delta^2}{2}}$ on $I_\delta$, from \eqref{eq:sjh} we can conclude that the distributional derivative of $h^\psi$ is 
zero in $I_\delta$. This means that there exists a real number $\beta(\psi)$ such that
\begin{equation}\label{eq:fond}
h^\psi(s)=\beta(\psi),\qquad \text{for }\LL^1-\text{a.e.~} s\in I_{\delta}.
\end{equation} 

Let $\Omega^\prime\subset\subset\Omega$, and select  
$\{\psi_i\} \subset C^1_c(\Omega\times \mathbb S^2)$ 
such that $\{\psi_i\}$ is dense in $C^0(\overline{\Omega^\prime}\times 
\mathbb S^2)$. Fix $\psi_i$, and
choose  $\eta_\delta\in C^\infty_c([-1,1])$ such that $0\leq\eta_\delta\leq 1$ on $[-1,1]$, $\eta_\delta\equiv 1$ on $I_{\delta/2}$.  
Before proceeding 
any further, let us recall that, by \cite[Proposition 3.4]{RS} (see also \cite[Lemma 4.4]{Chen}) there exists $\delta_0>0$  
independent  of $\eps$, such that if $\delta\leq\delta_0$ 
\begin{equation*}
\muepsueps(\Omega\cap\{\vert u_\eps\vert>1-\delta\})\leq C\delta,
\end{equation*}
where $C$ depends on $\Omega^\prime$, but not on $\eps$.

We then have
\begin{align*}
&\int_{-1}^1\eta_\delta
\sqrt{2W}\beta(\psi_i)\,ds=\int_{-1}^1\eta_\delta \sqrt{2W}\lim_{\eps\to 
0^+} h_{\eps}^{\psi_i} \,ds
=\lim_{\eps\to 0^+}\int_{-1}^1\eta_\delta g_{\eps}^{\psi_i} \,ds
\\
=& \lim_{\eps\to 0^+}\int_{\Omega\cap\{\vert 
u_\eps\vert<1-\frac{\delta}{2}\}}\psi_i\sqrt{2W(u_\eps)}\vert\nabla u_\eps\vert\,dx 
\\
& + \int_{\Omega\cap\{1-\frac{\delta}{2}<\vert u_\eps\vert<1-\delta\}}\eta_\delta(u_\eps)\psi_i\sqrt{W(u_\eps)}\vert\nabla u_\eps\vert\,dx
\\
&+\int_{\Omega\cap \{\vert u_\eps\vert>1-\frac{\delta}{2}\}}\psi_i\eps\sqrt{2W(u_\eps)}\vert\nabla u_\eps\vert\,dx-\int_{\Omega\cap \{\vert u_\eps\vert>1-\frac{\delta}{2}\}}\psi_i\sqrt{W(u_\eps)}\vert\nabla u_\eps\vert\,dx
\\
=&\surftens
\int\psi_i(y,\xi)\,dV^0(y,\star\xi) +O(\delta)
\\
=&\int_{-1}^1\sqrt{2W}\,ds V^0(\psi_i)+O(\delta)
=\int_{-1}^1\eta_\delta\sqrt{2W}V^0(\psi_i)\,ds
\\
&+\left( \int_{-1}^{-1+\frac{\delta}{2}}(1-\eta_\delta)\sqrt{2W}\,ds
   +  \int_{1-\frac{\delta}{2}}^1(1-\eta_\delta)\sqrt{2W}\,ds\right)
V^0(\psi_i)+O(\delta)
\\
=&\int_{-1}^1\eta_\delta\sqrt{2W}V^0(\psi_i)\,ds+O(\delta).
\end{align*}
Sending $\delta\to 0^+$ we obtain
\begin{equation}\label{eq:dragomerdo}
\int_{-1}^1\sqrt{2W}\,ds
~ \beta(\psi_i)= 
\int_{-1}^1\sqrt{2W}\,ds
~ V^0(\psi_i).
\end{equation}

Repeating the same argument for every $\psi_i$, by the density 
of $\{\psi_i\}$ in $C^0(\overline{\Omega^\prime}\times 
\mathbb S^2)$ and \eqref{eq:dragomerdo} we deduce that $ \beta=V^0$ as measures on $G^0_2(\Omega^\prime)$. 
\endproof
Let $\psi\in C^1_c(\Omega\times \mathbb S^2)$.
From the estimates on $(d/ds)g^\psi_\eps$ obtained in the proof of Lemma \ref{tuttoperCicciobelloRicchione} we can conclude that 
$g_\eps^\psi\to g^\psi$ strongly in $W^{1,1}_{{\rm loc}}((-1,1))$ as $\eps\to 
0^+$. The proof of Theorem \ref{theo:se2mni}-(B) is complete.

We are now in a position to conclude 
the proof of Theorem \ref{theo:se2mni}-(A).
\begin{proposition}\label{smandruppone}
There exists a 
(not relabelled) subsequence $\{\Vzeroepsueps\}$ converging,
as oriented varifolds, to 
$V^0=\var(\rectifiableset,
\tau,\theta_1,\theta_2)\in\mathbf{IV}^0_2(\Omega)$, with $\q_\sharp V^0=V$. 
\end{proposition}

\proof
As we already noticed at the beginning of the present section, by 
\eqref{eq:ener-bound}, we can extract a  subsequence such that 
$\Vzeroepsueps$ converge as Radon measures  to a Radon measure $V^0$ on 
$G^0_2(\Omega)$, and also that $\q_\sharp  V^0=V$. Hence, in order to conclude it remains to show that $V^0\in\mathbf{IV}^0_2(\Omega)$.
To this aim we will make use of Lemma \ref{tuttoperCicciobelloRicchione}.

Fix $\Omega^\prime\subset
\subset\Omega$ with smooth boundary.
\noindent By 
Sard's Lemma and Lemma \ref{tuttoperCicciobelloRicchione} we can find a  
subsequence $\{\Vzeroepsuepsk\}_{k}$ and  a subset $J\subset 
[-1,1]$, with $\LL^1(J)=0$, such that for every $s\in [-1,1]\setminus J$,
\begin{gather*}
\{u_{\eps_k}=s\}\text{ is a smooth embedded 
surface and } \{u_{\eps_k}=s\}\cap\{\nabla u_{\eps_k}=0\}=\emptyset,
\\
\partial\llbracket \var(\{u_{\eps_k}=s\},\star
\niuepsk
,1)\rrbracket(\Omega^\prime)=0,
\\
\lim_{k\to\infty} \var(\{u_{\eps_k}=s\},\star
\niuepsk
,1)=V^0~ \text{ as oriented varifolds on }\Omega^\prime.
\end{gather*}

Next we fix $\delta>0$ and set $I_\delta:=[-1+\delta,1-\delta]$.  Since
we have
\begin{align*}
&\int_{I_\delta\setminus J}\left\vert\delta \var(\{u_{\eps_k}=s\},\star
\niuepsk
,1)\right\vert(\Omega^\prime)\,ds=
\int_{I_\delta\setminus J}\int_{\{u_{\eps_k}=s\}\cap\Omega^\prime}\left\vert \mathrm{div}\left(
\niuepsk
\right)\right\vert\,d\Ha^2ds
\\
\leq & \frac{1}{(2\delta-\delta^2)}\int_{\Omega^\prime}\left\vert \mathrm{div}\left(
\niuepsk
\right)\right\vert\sqrt{2W(u_{\eps_k})}\vert\nabla u_{\eps_k}\vert\,dx\leq\frac{2}{(2\delta-\delta^2)}
\int_{\Omega^\prime}\vert\sff_{u_{\eps_k}}\vert\sqrt{2W(u_{\eps_k})}\vert\nabla u_{\eps_k}\vert\,dx
\\
\leq &\frac{2}{(2\delta-\delta^2)} \left(\int_{\Omega}\vert\sff_{u_{\eps_k}}\vert^2\;d\tildemuepsuepsk\right)^{1/2}
\left\{\left[
\tildemuepsuepsk(\Omega)\right]^{1/2}+
2\left[\vert \xiepsuepsk\vert(\Omega)\right]^{1/2}\right\},
\end{align*}
by the choice 
of the $\eps_k$, the set $J$ and \eqref{eq:granmatthias}, we can conclude that there exists 
$s = s_{\eps_k}\in I_\delta\setminus J$ such that
\begin{equation*}
\limsup_{k\to \infty}\Big\vert\delta \var(\{u_{\eps_k}=s_{\eps_k}\},\star
\niuepsk
,1)\Big\vert(\Omega^\prime)<+\infty.
\end{equation*}
The thesis is then a direct consequence of the properties of 
$\{u_{\eps_k}=s_{\eps_k}\}$ for $s\in I_\delta\setminus J$ and Theorem 
\ref{the:comp-or-var}. 
\endproof
%
\section{Proof of Proposition \ref{lem:conv-discr-BV}}\label{WUgoCenci}  
As in Section \ref{maanvediunpodanna}, by \eqref{eq:ener-bound} 
we  deduce that \eqref{eq:granmatthias} holds. Hence we can apply Theorem \ref{the:RS} and conclude 
that, up to selecting a further subsequence, \eqref{eq:the-RS-van-discr} 
holds. In addition, the densities of the discrepany measures
are uniformly bounded in $L^1(\Omega)$, and we have 
\begin{equation*}
\begin{split}
&\int_\Omega\left\vert\nabla 
\xi_{u_\eps}^\eps
\right\vert\,dx=
\int_\Omega\left\vert\eps\nabla^2u_\eps\nabla 
u_\eps-\frac{W^\prime(u_\eps)}{\eps}
\nabla u_\eps\right\vert\,dx
\\
=&\int_{\{\nabla u_\eps\neq 0\}}\left\vert\left[\eps\nabla^2 u_\eps
-\frac{W^\prime(u_\eps)}{\eps}\frac{\nabla u_\eps\otimes\nabla 
u_\eps}{\vert\nabla u_\eps\vert^2}\right]\nabla u_\eps\right\vert\, dx
\\
\leq&3^{1/4}\left(\frac{1}{\eps}\int_{\{\nabla u_\eps\neq 
0\}}\left\vert\eps\nabla^2u_\eps
-\frac{W^\prime(u_\eps)}{\eps}\frac{\nabla u_\eps\otimes\nabla u_\eps}{\vert\nabla u_\eps\vert^2}\right\vert^2\, dx\right)^{1/2}
\left(
\tildemuepsueps(\Omega)\right)^{1/2}
\\
=&3^{1/4}\left(\int_\Omega\vert B^\eps_{u_\eps} 
\vert^2\,d
\tildemuepsueps
\right)^{1/2}
\left[\tildemuepsueps
(\Omega)\right]^{1/2}\leq \constant,
\end{split}
\end{equation*}
where $\constant$ is a positive constant independent of $\eps$.

By the compactness theorem in $BV$ (see \cite{AFPa00}) and Theorem \ref{the:RS}  we can select a 
further  subsequence  such that $\xi_{u_\eps}^\eps
\rightharpoonup 0$ 
weakly in $BV(\Omega)$ as $\eps\to 0^+$. Moreover \eqref{asburatto} holds 
by Rellich-Kondrachov compactness theorem (see \cite{AFPa00}).
\endproof
%

\section{Final Comments}\label{bombemerda}
\subsection{Relaxing the constraints on $\elbend,\,\elGauss$} \label{Del-puzzone}
As already stated in Remark \ref{rem:Sil-barb}, 
Theorem \ref{the:Helf-lsc} still holds when replacing \eqref{eq:constr} 
with the more general constraint $-2<\elGauss/\elbend<0$. 
Although we cannot prove 
Theorem \ref{the:Gamma-liminf-cond} 
(and hence Corollary \ref{cor:Gamma-conv-smooth-points})
when $-2<\elGauss/\elbend<0$, 
we can relax 
condition \eqref{eq:constr} 
to 
\begin{equation}\label{ore}
\elGauss<0<\frac{3}{2}\elbend+\elGauss.
\end{equation}
In fact, in this case we can still derive  \eqref{eq:granmatthias} 
using the inequality
\begin{equation*}
(\fepsu)^2
=
\left(\mathrm{tr}(\calBepsu)\right)^2
\leq 3\vert\calBepsu\vert^2.
\end{equation*}
 Hence, in particular, Theorem \ref{the:Gamma-liminf-cond} holds for  
$
\elbend=-\elGauss=1,
$
which gives the usual isotropic bending energy 
\begin{gather*}
\Helf(E)=\frac{1}{2}
\int_{\Omega \cap \partial E}
\vert\sff_{\partial E}\vert^2\,d\Ha^2,
\\
\mathpzc W_\eps(u) =
\frac{1}{2\eps}\int_\Omega\left\vert\eps\nabla^2u-
\frac{W^\prime(u)}{\eps}\niu\otimes\niu\right\vert^2\,dx.
\end{gather*}
%
\subsection{Full $\Gamma$-convergence and convergence of constrained minimizers.}\label{alvaro}
Corollary \ref{cor:Gamma-conv-smooth-points} shows that the $\Gamma$-limit with respect to the $L^1$-topology of $\mathpzc W_\eps$ is given by $\Helf$ on smooth points.
However, since $\Gamma$-limits are always lower semicontinuous, the natural candidate 
for a full $\Gamma$-convergence result is 
the $L^1$-lower semicontinuous envelope $\overline{\Helf}$ 
of $\Helf$ 
defined by
\begin{multline*}
\overline{\Helf}(E):=\inf
\Big\{\liminf_{h\to\infty}\Helf(E_h):\, E_h \subset \Omega 
~{\rm bounded~with~} \partial E_h \in 
C^2,\,
\\
 \lim_{h\to\infty}\chi_{E_h}
=\chi_E \text{ in }L^1(\Omega)\Big\}.
\end{multline*}
Let us recall some facts about $\overline{\Helf}$ (see for example \cite{Del}). Define
\begin{multline*}
\mathcal D:=\Big\{W\in\mathbf{IV}_2(\Omega): \,W=\lim_{h\to\infty}
\var(\partial E_h,1),\, 
E_h \subset \Omega {\rm ~bounded~with~}
\partial E_h \in C^2,\,
\\
\sup_{h\in\N}\int_{\Omega \cap 
\partial E_h}
\left[1+\vert\sff_{\partial E_h}\vert^2\right]\,d\Ha^2<+\infty\Big\},
\end{multline*}
and
\begin{multline*}
 \mathcal A(E):=\Big\{ W\in \mathcal D: \,W=\lim_{h\to\infty}
\var(\partial E_h,1),\, 
E_h \subset \Omega {\rm ~bounded~with~}
\partial E_h \in C^2,\,
\\
\lim_{h\to\infty}\chi_{E_h}=\chi_E\text{ in }L^1(\Omega)\Big\}.
\end{multline*}
Eventually, we recall that if 
$W\in\mathcal D$ then $W\in\mathcal A(E_{W})$ where $E_{W}$ is an open, bounded subset with finite perimeter in $\Omega$, such that the essential boundary of $E$ coincides with the set of points of odd $2$-density with respect to $\mu_{W}$.

{}From \cite[Corollary 5.4]{Del}, we obtain
\begin{equation*}
\overline{\Helf}(E)=\min\{\Helf(V):V\in\mathcal A(E)\}.
\end{equation*}
Hence, if we would be able to prove that $V=\lim_{\eps\to 0^+} 
\Vepsueps\in \mathcal A(E)$, by \eqref{eq:Gamma-liminf-var} we would have
\begin{equation*}
\liminf_{\eps \to 0^+}
\mathpzc W_\eps(u_\eps)\geq \surftens
 \Helf (V)\geq\surftens
\overline{\Helf}(E),
\end{equation*}
which, together with $\overline{\Helf}(E)=\Helf(E)$ for 
$E\subset \Omega$ bounded with boundary of class $C^2$,
 would imply that $\Gamma(L^1(\Omega))-\lim_{\eps\to 0^+}\mathpzc W_\eps=\overline{\Helf}$.  Although Theorem \ref{theo:se2mni}-(B) seems to represent a signicative step in this direction, in order to prove that $V\in\mathcal A(E)$ we miss an estimate similar to the one
 proved in \cite[Lemma 2]{Tonesff}, \cite[Theorem 1]{ToneNotes}.
Actually, we are able to prove that $V\in\mathcal A(E)$ 
under the stronger assumption
\begin{gather}
\sup_{0<\eps<1}\widetilde{\mathpzc W}_\eps(u_\eps)<+\infty,
\label{eq:padre-peo}
\\
\widetilde{\mathpzc W}_\eps(u_\eps):=\mathpzc W_\eps(u_\eps)
+\int_\Omega \vert B^\eps_{u_\eps}\vert^2\frac{W(u_\eps)}{\eps}\,dx.
\notag
\end{gather}
Indeed, assuming that \eqref{eq:padre-peo} holds, we have
\begin{gather*}
\sup_{0<\eps<1}\int_{-1}^1
\int_{\{u_\eps=s\}}\vert\sff_{\{u_\eps=s\}}\vert^2\,d\Ha^2ds\leq
\surftens^{-1}
\sup_{0<\eps<1}\widetilde{\mathpzc W}_\eps(u_\eps)<+\infty,
\end{gather*}
which, by Lemma \ref{tuttoperCicciobelloRicchione}, gives $V\in\mathcal A(E)$. Moreover, this means that we can conclude that 
chosen $\widetilde u_\eps$ so that
\begin{gather*}
\widetilde{\mathpzc W}_\eps(\widetilde u_\eps)=
\min\left\{\widetilde{\mathpzc W}_\eps(u):\, \mathpzc P_\eps(u)=\Lambda_1,\, \int_{\Omega}\frac{1+u}{2}\,dx=\Lambda_2\right\}
\end{gather*}
we have, up to a subsequence,
\begin{equation*}
V^\eps_{\widetilde u_\eps}\to \widetilde{V}\in\mathcal D,
\qquad u_\eps\to u=2\chi_{\widetilde E}-1 \qquad {\rm as}~ \eps \to 0^+,
\end{equation*}
where 
\begin{itemize}
\item[-] $\widetilde{V}$ solves
\begin{gather*}
\min\Big\{\Helf(V):\,V\in\mathcal D, ~ \mu_{V}(\Omega)=\Lambda_1, ~ \mathcal 
L^3(\Omega \cap E_{V})=\Lambda_2\Big\}
\end{gather*}
\item[-] $\widetilde E\subset\Omega$ solves
\begin{multline*}
\min\left\{\overline{\Helf}(E):\,
  \forall  
W\in\mathcal A(E)\text{ we have } \mu_{W}(\Omega)
=\Lambda_1,\, \mathcal L^3(\Omega\cap E)=\Lambda_2\right\}.
\end{multline*}
\item[-]  $\Helf(\widetilde V)=\overline{\Helf}(\widetilde E)$.
\end{itemize}
\subsection{The case of non-zero spontaneous curvature.}\label{ilcamionista}
As we already remarked in the introduction, when $H_0\neq 0$ the functional
\begin{gather}\label{eq:cagnone}
\int_{\partial E\cap \Omega}(H_{\partial E}-H_0)^2 ~d\Ha^2
\end{gather}
not only depends  on the surface $\partial E$ but also on the orientation of $\partial E$. Moreover such a 
functional is not lower semicontinuous with respect to the varifolds convergence. In fact, as an example due to Karsten Gro{\ss}e-Brauckmann shows (see \cite{Gross}, \cite{GrossPhd} and \cite{SchaJDG}), 
there exists
a sequence $\{E_h\}_h$ 
of 
smooth sets 
in $\Omega:=B(0,1)$, such that for every $h\in\N$ the surface $\partial E_h$ has constant (scalar) mean curvature equal to $1$, and at the same time the sequence of  varifolds $\var(\partial E_h,1)$ converges to the 
varifold $\var(
\langle\mathbf{e}_3
\rangle^\perp,2)
$ 
 in $\Omega$. Hence, assuming $H_0=1$, we have
$$
0=\lim_{h\to\infty}\int_{\Omega \cap \partial E_h}
(H_{\partial E_h}-H_0)^2\,d\Ha^2<2\pi=2\int_{\
\langle\mathbf{e}_3
\rangle^\perp
\cap B(0,1)}(H_0)^2
~d\Ha^2.
$$
However if we consider the complete Helfrich's energy
\begin{equation}\label{eq:cagnone-bastardo}
\Helf(E)=\int_{\Omega \cap 
\partial E}
\left[
\frac{\elbend}{2}\Big(H_{\partial E}-H_0\Big)^2+\elGauss K_{\partial E}
\right]
~d\Ha^2,
\end{equation}
and assume (as in the case of zero spontaneous curvature) that $-2<\elbend/\elGauss<0$, the results of \cite{Del} still apply and Theorem \ref{the:Helf-lsc} holds also in this case. Moreover the functional is lower
semicontinuous with repect to the convergence of the oriented varifolds and, whenever $\sup_{h\in\N}\Helf(E_h)<+\infty$, the oriented varifolds $\var(\partial E_h,\,\star\nu_{\partial E_h},\,1)$ converge (up to a subsequence) to an \textit{oriented curvature varifold} $V^0\in\mathbf{IV}^0_2(\Omega)$ in the sense of \cite{DelOriVar}. 

Possible diffuse-interface approximating functionals
 for \eqref{eq:cagnone} are 
\begin{gather}\label{eq:daddario}
\frac{1}{\eps}\int_\Omega \Big(f^\eps_u-H_0\eps\vert\nabla u\vert\Big)^2\,dx,
\quad
\frac{1}{\eps}\int_\Omega\Big(
f_{u}^\eps-H_0\sqrt{2W(u)}\Big)^2\,dx,
\end{gather}
the latter being the one proposed in \cite{DuSpontCurv}. Consequently a natural candidate for the diffuse-interface approximation of \eqref{eq:cagnone-bastardo}
is
\begin{gather*}
\widehat{\mathpzc W}_\eps(u):=\frac{\elbend}{2}\widehat{\mathpzc H}_\eps(u)+\elGauss\mathpzc K_\eps(u),
\end{gather*}
where $\widehat{\mathpzc H}_\eps(u)$ is given by one of the two expressions in \eqref{eq:daddario}.
If \eqref{eq:constr} is satisfied, by a
direct calculation we can show that \eqref{eq:granmatthias} holds as soon as
$$
\sup_{0<\eps<1}\Big(\mu^\eps_{u_\eps}(\Omega)+\widehat{\mathpzc W}_\eps(u_\eps)\Big)<+\infty.
$$
Hence we can conclude that also Lemma \ref{lem:Kazumba} and Lemma \ref{lem:conv-to-Hutch} apply and, with minor modifications to the arguments of Sections \ref{maanvediunpodanna}-\ref{WUgoCenci},
we can prove that Theorem \ref{theo:se2mni} and Proposition \ref{lem:conv-discr-BV}
hold also for $\widehat{\mathpzc W}_\eps$. Moreover we can use the same sequence  $\{\widehat u_\eps\}_\eps\subset C^2(\Omega)$ constructed in Section \ref{prendodimiraquestalberoinmezzoalgiardino} to show that also an 
analog  of  Theorem \ref{the:Gamma-limsup} holds 
for $\widehat{\mathpzc W}_\eps$. However, in order to prove that the lower bound estimate corresponding to \eqref{eq:Gamma-liminf-var} holds, 
we should prove that 
\begin{gather}\label{eq:Ossesionatore-Romano}
\lim_{\eps\to 0^+}\int_\Omega\nabla\xi^\eps_{u_\eps}\cdot\nu_{u_\eps}\,dx=0.
\end{gather}
Unfortunately we are not able to prove \eqref{eq:Ossesionatore-Romano} 
unless additional hypothesis are made on $u_\eps$ 
(for example if $\muepsueps\to 2\surftens\vert\nabla\chi_E\vert$, then \eqref{eq:Ossesionatore-Romano} follows from \eqref{Kazumba}, Theorem 
\ref{theo:se2mni} and Lemma \ref{lem:Moser}). 
 However, a possible strategy to obtain \eqref{eq:Ossesionatore-Romano} might 
be  trying to use  Proposition \ref{lem:conv-discr-BV} on each of the ``well-separated transition layers'' that can be 
obtained via an appropriate blow-up procedure  (see  \cite[Proposition 5.3]{RS}), and then conclude via a covering argument.
%
\begin{appendix}
\section{Measure-function pairs}\label{trappolapernani}

Let $D\subset\R^l$; 
we say that $(\mu,f)$ is a \textit{measure-function pair over $D$ 
with values in $\R^m$}, if 
 $\mu$ is a
positive Radon measure on 
$D$, 
$f:D\to\R^m$ is defined $\mu$-almost
everywhere and $f\in L^1_{{\rm loc}}(\mu)$. 

Let us 
 recall the definition of measure-function pairs convergence (see \cite{Hu})
\begin{definition}
Let $(\mu_k,f_k)$, $(\mu,f)$ be measure-function pairs
 on $D$ with values in $\R^m$ for every $k\in\N$. We say that 
$(\mu_k,f_k)$ converge weakly 
to $(\mu,f)$ as measure-function pairs as $k\to\infty$ if 
\begin{equation*}
\lim_{k\to\infty}\int f_k\cdot Y\,d\mu_k=\int f\cdot Y\, d\mu\qquad \forall Y\in C^0_c(D,\R^m).
\end{equation*}
\end{definition}
\begin{definition}
We say that a function $F:\R^m\to[0,+\infty)$ is a standard integrand 
provided 
$F$ is strictly convex on $\R^m$, and 
\begin{equation*}
g(\vert q\vert)\leq F(q)\qquad \forall q\in\R^m,
\end{equation*}
where 
$g\in C^0([0,+\infty))$ is non-negative,  
 increasing and $g(t)\to+\infty$ as $t\to+\infty$.
\end{definition}
\begin{definition}\label{def:strong-meas-func-conv}
Let $(\mu_k,f_k)$ and $(\mu,f)$ 
be measure-function pairs over $D$ with values on $\R^m$. Suppose $\mu_k
\rightharpoonup\mu$ as $k\to\infty$ as 
Radon measures. We say that $(\mu_k,f_k)$ converge to $(\mu,f)$ 
in the $F$-strong sense in $D$ if 
\begin{itemize}
\item[(i)] $\int F(f_k)\,d\mu_k<+\infty$ for every $k\in\N$;
\item[(ii)] setting $D_{kj}:=\{y\in D:\, \vert f_k(y)\vert\geq j\}$ we have
$$
\lim_{k\to\infty}\int_{D_{kj}}F(f_k)\,d\mu_k=0,
$$
uniformly in $k\in\N$;
\item[(iii)] for every $\psi\in C^0_c(D \times \R^m)$ we have
$$
\lim_{k\to\infty}\int \psi(y,f_k)\,d\mu_k=\int \psi(y,f)\,d\mu.
$$
\end{itemize}
\end{definition}
We say that a sequence of measure-function pairs converges $L^p$-strongly 
($p\in [1,\infty)$) if it converges strongly in the $F_p$-sense, 
with $F_p(q):=\vert q\vert^p$.

The following result has been proved in \cite[Theorem 4.4.2]{Hu}.
\begin{theorem}\label{the:Hutch-comp-lsc}
Let $(\mu_k,f_k)_{k\in\N}$ be measure-function pairs over $D$ 
with values in $\R^m$. Suppose that $\mu$ is a Radon measure on $D$ 
and $\mu_k\rightharpoonup\mu$ in $D$ 
as $k\to\infty$. Let $F:\R^m\to[0,+\infty)$ be a standard integrand. 
The following assertions hold.
\begin{itemize}
\item[(i)] If
\begin{equation}\label{eq:miao}
\sup_{k\in\N}\int F(f_k)\,d\mu_k<+\infty,
\end{equation}
then there exists $f\in L^1_{{\rm loc}}(\mu)$ and a (not
relabelled) subsequence $\{(\mu_k,f_k)\}$ such that 
\begin{equation}\label{eq:bau}
\lim_{k\to\infty}(\mu_k,f_k)=(\mu,f),
\end{equation}
weakly as measure-function pairs on $D$ with values on $\R^m$.
\item[(ii)] If $\{(\mu_k,f_k)\}$ and $(\mu,f)$ satisfy \eqref{eq:miao}, \eqref{eq:bau}, 
then
\begin{equation}\label{eq:lsc-the-Hutch}
\int F(f)\, d\mu\leq \liminf_{k\to\infty}\int F(f_k)\,d\mu_k.
\end{equation}
\end{itemize}
\end{theorem}

\begin{remark}\rm
We can adapt the notions and results proved until this point in the present Appendix  to the case where $D$ is an open subset of a smooth manifold embedded in $\R^m$ for some $m\in\N$. In particular, in our applications we will often consider $D=G_2(\Omega)$ or $D=G^0_2(\Omega)$.
\end{remark} 
The following lemma is a particular case of \cite[Proposition 3.2]{Mose01}.
\begin{lemma}\label{lem:Moser}
Let $(\mu_k,g_k)$ and 
$(\mu, g)$ be measure-function pairs on $D$ with values in $\R^m$ such that 
$$
\sup_{k\in\N}\|g_k\|_{L^2(\mu_k)}<+\infty,
$$ 
and $(\mu_k,g_k)$ weakly converge to $(\mu,g)$ as measure-function pairs. 

Moreover let $(\mu_k,f_k),\,(\mu,f)$ be 
measure-function pairs on $D$ with values in $R^m$  such that $(\mu_k,f_k)$ converges $L^2$-strongly to $(\mu,f)$. Then 
$$
\lim_{k\to\infty}(\mu_k, f_k\cdot g_k)=(\mu,f\cdot g),
$$
weakly as measure-function pairs on $D$ with values in $\R$.
\end{lemma}
\section{Geometric Measure Theory: varifolds}\label{luponeassassino}
Let us recall some basic fact in the theory of 
varifolds, the main bibliographic sources being \cite{Si} and \cite{Hu}.

We call \textit{varifold} (resp. \textit{oriented varifold}) any positive 
Radon measure on $G_2(\Omega)$ (resp. on $G^0_2(\Omega)$). 
In this paper we are confined to surfaces, hence 
we use the terms varifold and oriented varifold to mean a $2$-varifold in $\Omega$.

If $V^0$ is an oriented 
varifold then the push-forward
$\q_\sharp V^0$
is the corresponding unoriented 
varifold associated with $V^0$ by projection onto $G_2(\Omega)$.

For any varifold (or oriented varifold) $V$ we define 
$\mu_V$ to be the Radon  measure on $\Omega$  
obtained by projecting $V$ onto $\Omega$.

Let $\rectifiableset
$ be a $2$-rectifiable subset of $\R^3$ with  finite $\Ha^2$-measure
and let  $\theta,\,\theta_1,\,\theta_2:\rectifiableset
\to\R^+$ be 
$\Ha^2 \res \rectifiableset
$-measurable 
functions.  Suppose $\tau:\rectifiableset
\to G^0_{2,3}$ 
is $\Ha^2\res \rectifiableset
$-measurable 
and $\q(\tau(x))=T_x \rectifiableset
$ for $\Ha^2\res 
\rectifiableset
$-almost everywhere $x$ ($\tau$ is called an orientation function on $\mathcal 
M$). Then we define the \textit{rectifiable}  
(unoriented and oriented respectively) \textit{varifolds}
$$
V=\var(\rectifiableset
,\theta),\qquad V^0=\var(\rectifiableset
,\tau,\theta_1)
+\var(\rectifiableset
,-\tau,\theta_2) =: 
\var(\rectifiableset, \tau, \theta_1,\theta_2),
$$
by
\begin{gather*}
V(\phi):=\int_{\rectifiableset
}
\phi(x,T_x \rectifiableset
)\,\theta(x)d\Ha^2 \qquad\forall\phi\in 
C^0_c(G_2(\Omega)),
\\
V^0(\varphi):=\int_{\rectifiableset
}
\Big[\varphi(x,\tau(x))\theta_1(x)+\varphi(x,-\tau(x))\theta_2(x) \Big]\,
d\Ha^2\qquad \forall\varphi\in C^0_c(G^0_2(\Omega)).
\end{gather*}
With the notation
$\var(\rectifiableset, \tau, \theta)$ we mean 
$\var(\rectifiableset, \tau, \theta,0)$.

When $\theta$ (resp. $\theta_1,\,\theta_2$) take values in $\N$
 we say that $V=\var(\rectifiableset
,\theta)$ (resp. $V^0=\var(\rectifiableset, \tau, 
\theta_1,\theta_2)$) is a  \textit{rectifiable integer}  
unoriented (resp. oriented) \textit{varifold}
and we write $V \in \mathbf{IV}_2(\Omega)$ (resp. 
$V^0 \in \mathbf{IV}^0_2(\Omega)$).
If $V^0 = {\bf v}(\mathcal M, \tau, 
\theta_1,\theta_2)\in\mathbf{IV}^0_2(\Omega)$ 
the integral rectifiable $2$-current $\llbracket V^0\rrbracket$ is defined as
\begin{equation*}
\llbracket V^0\rrbracket(\omega):=
\int_M \langle \omega(x),\tau(x)\rangle\, \left(\theta_1(x)-\theta_2(x)\right)d\Ha^2(x)
\qquad \forall \omega\in C^0(\Omega,\Lambda_2(\R^3)). 
\end{equation*}
As usual $\partial\llbracket V^0\rrbracket$ denotes
 the boundary of the current $\llbracket V^0\rrbracket$,
and $\vert \partial \llbracket V^0\rrbracket\vert$
is the mass of $\llbracket \partial V^0\rrbracket$ (see \cite{Si}).

Let $V$ be an unoriented 
varifold on $\Omega$; we define \textit{the first variation 
of} $V$ as the linear operator
\begin{equation*}
\delta V:C^1_c(\Omega,\R^3)\to\R,\qquad Y\to \int {\rm tr}(S\nabla 
Y(x))\,dV(x,S).
\end{equation*}

We say that $V$ has \textit{bounded first variation}  
(resp. 
\textit{generalized mean curvature in} $L^p$, $p>1$) if $\delta V$ 
can be 
extended to a linear continuous operator on  $C^0_c(\Omega,\R^3)$ 
(resp. on $L^p(\mu_V,\R^3)$). In this case $\vert \delta V\vert$ 
denotes the total variation of $\delta V$.
Whenever the 
varifold $V$ has bounded first variation we  call \textit{the generalized 
mean curvature vector of} $V$ the vector field 
$$
\Hv_{V}=\frac{d\delta 
V}{d\mu_V},
$$
where the right-hand side denotes the Radon-Nikodym derivative.

By varifold convergence (resp. oriented varifold convergence) we 
mean the convergence as Radon measures
 on $G_2(\Omega)$ (resp. on $G^0_2(\Omega)$). 
The following compactness theorem for oriented varifolds 
is proved in \cite[Theorem 3.1]{Hu}.
\begin{theorem}\label{the:comp-or-var}
Let $\constant>0$ and let $\{\Omega_i\}$ be a sequence of open subsets 
with smooth boundary invading $\Omega$. The set
\begin{equation*}
\left\{V^{0}\in \mathbf{IV}_2^0(\Omega):\,\forall i\in\N,\, \mu_{\q_\sharp(V^{0})}(\Omega_i)+
\vert\delta 
(\q_\sharp V^{0})\vert(\Omega_i)+\vert\partial \llbracket V^{0}\rrbracket\vert(\Omega_i)\leq 
\constant \right\}
\end{equation*}
is sequentially compact with respect to the oriented varifolds convergence.
\end{theorem}
\begin{remark}\label{rem:lerchiopattino}
Let $\{V_h\}$ be a sequence of varifolds converging to a varifold $V$,
and suppose that there exist $\mu_{V_h}$-measurable maps 
$S^h_\cdot$ 
and a $\mu_V$-measurable map 
$S_\cdot$ such that
\begin{align*}
&V_h(\Psi)=\int\Psi(x,S^h_x)\,d\mu_{V_h}(x)
\qquad \forall\Psi\in C^0_c(G_2(\Omega)),~\forall h\in\N
\\
&V(\Psi)=\int\Psi(x,S_x)\,d\mu_{V}(x) \qquad \forall\Psi\in C^0_c(G_2(\Omega)).
\end{align*}
Then it can be 
checked that the measure function pair $(\mu_{V_h},S^h_\cdot)$ converge $L^p$-strongly to $(\mu_V, S_\cdot)$ as measure function pairs on $\Omega$ 
with values in $G_2(\Omega)$, for every $p\in (1,+\infty)$.
\end{remark}
Following \cite{Hu} we define the notion of Hutchinson's
curvature varifold with generalized second fundamental form.

\begin{definition}\label{def:Hutch-var}
Let $V \in {\bf IV}_2(\Omega)$.
We say that 
$V$ is a curvature varifold with generalized second 
fundamental form in  $L^2$, if there exists $A^V=A^V_{ijk}\in 
L^2(V,\R^{3^3})$  such that for every 
function $\phi\in C_c^1(G_2(\Omega))$ and $i=1,2,3$, 
\begin{equation}\label{eq:def-Hutch-var} 
\int_{G_{2}(\Omega)}(S_{ij}\partial_j 
\phi+A^V_{ijk}D_{m_{jk}}\phi+A^V_{jij}\phi)\, 
dV(x,S)=0,
\end{equation}
where $D_{m^{}_{jk}} \phi$ denotes the derivative of $\phi (x,\cdot)$ 
with respect to its $jk$-entry variable.

Moreover we define the generalized second fundamental form $\sff_V
=(B^k_{ij})_{1\leq i,j,k\leq 3}$ of $V$ as
\begin{equation}\label{eq:def-gen-sff}
B^k_{ij}(x,S):=S_{jl}A^V_{ikl}(x,S).
\end{equation}
\end{definition}
\begin{remark}\label{rem:MC-Hutch-MC-Allard}\rm
Every curvature varifold $V$ with generalized  second
fundamental form in $L^2$ has bounded first variation. Moreover
\begin{equation}\label{eq:MC-Hutch-MC-Allard} 
\Hv_V(x)=(A_{j1j}(x,T_x\mu_V),A_{j2j}(x,T_x\mu_V),A_{j3j}(x,T_x\mu_V))\in L^2(\mu_V,\R^3),
\end{equation}
for $\mu_V$ almost every $x \in \Omega$.
\end{remark}
\begin{remark}\rm
If $V=\var(\Sigma,1)$, where $\Sigma$ is a smooth, compact 
surface without boundary, the generalized second fundamental form as well as the mean curvature and the tensor $A_V$ coincide with the classical quantities defined in Section \ref{prel:diff-geom}, and the same is true for the oriented varifold associated with
 $\Sigma$. Moreover the generalized second fundamental form
 and the functions $A^V_{ijk}$ verify Proposition \ref{prop:A-vs-B}. 
\end{remark} 
Next we give a definition of convergence for Hutchinson's
curvature varifolds.
\begin{definition}\label{def:Hutch-conv}
Let $\{V_h\}$ be a sequence of 
curvature varifolds with  generalized 
second fundamental form in $L^2$, and let $V$ 
be a curvature varifold with  generalized
second fundamental form in $L^2$. 
We say that $V_h$ converge as curvature varifolds to 
$V$  if
\begin{align*}
& \lim_{h\to\infty}V_h=V \qquad{\it ~as~ varifolds},
\\
& \lim_{h\to\infty}(V_h,A_{V_h})=(V,A_V)\qquad {\it ~as ~measure-function 
~pairs}.
\end{align*}
\end{definition}
\begin{remark}\rm\label{rem:conv-A-conv-sff}
By Remark \ref{rem:lerchiopattino}, Lemma \ref{lem:Moser} and the definition of generalized second fundamental form $\sff_{V_h}$, we have that if $V_h\to V$ as curvature varifolds then 
\begin{equation*}
(V_h,\sff_{V_h})\to (V,\sff_V)
\end{equation*}
as measure-function pairs on $G_{2}(\Omega)$ with values in $\R^{3^3}$.
\end{remark} 

As a consequence of Definition \ref{def:Hutch-conv} and Theorem 
\ref{the:Hutch-comp-lsc} we have the following
\begin{proposition}\label{prop:compactness-curv-var} 
Let $\{V_h\} \subset\mathbf{IV}_2(\Omega)$ 
 be a sequence of  
curvature varifolds 
with generalized second fundamental form in $L^2$
satisfying 
\begin{equation*}
\begin{split}
&~~\sup_{h\in\N}\left\{\mu_{V_h}(\Omega)+\int \sum_{i,j,k=1}^3
(A^{V_h}_{ijk})^2\,dV_h<+\infty\right\}.
\end{split}
\end{equation*}
Then $\{V_h\}$ has a subsequence converging to $V\in\mathbf{IV}_2(\Omega)$ 
as curvature varifolds.
\end{proposition}
\end{appendix}
\bibliography{Bel_Mug_09}

\begin{thebibliography}{10}

\bibitem{AFPa00}
Luigi Ambrosio, Nicola Fusco, and Diego Pallara.
\newblock {\em Functions of bounded variation and free discontinuity problems}.
\newblock Oxford Mathematical Monographs. The Clarendon Press Oxford University
  Press, New York, 2000.

\bibitem{ArrDeSim}
Marino Arroyo and Antonio De~Simone.
\newblock Relaxation dynamics of uid membranes.
\newblock {\em Phys. Rev. E}, 79(3):0319151--03191517, 2009.

\bibitem{BaumgartNature}
Tobias Baumgart, Samuel~T. Hess, and Webb~E. Webb.
\newblock Imaging coexisting fluid domains in biomembrane models coupling
  curvature and line tension.
\newblock {\em Nature}, 425:821--824, 2003.

\bibitem{BMRad}
Giovanni Bellettini and Luca Mugnai.
\newblock On the approximation of the elastica functional in radial symmetry.
\newblock {\em Calc. Var. Partial Differential Equations}, 24(1):1--20, 2005.

\bibitem{BePa}
Giovanni Bellettini and Maurizio Paolini.
\newblock Approssimazione variazionale di funzionali con curvatura.
\newblock {\em Seminario Analisi Matematica Univ. Bologna, Tecnoprint}, pages
  87--97, 1993.

\bibitem{AmAll}
Martine Ben~Amar and Jean~Marc Allain.
\newblock Budding and fission of a multiphase vesicle.
\newblock {\em Eur. Phys. J. E}, 20:409--420, 2006.

\bibitem{BiKaMi}
Thierry Biben, Klaus Kassner, and Chaouqi Misbah.
\newblock Phase-field approach to three-dimensional vesicle dynamics.
\newblock {\em Physical Review E}, 72:041921, 2005.

\bibitem{Boal}
David Boal.
\newblock {\em Mechanics of the Cell}.
\newblock Cambridge University Press, Cambridge, 2002.

\bibitem{Braides-book}
Andrea Braides.
\newblock {\em $\Gamma$-convergence for beginners}, volume~22 of {\em Oxford
  Lecture Series in Mathematics and its Applications}.
\newblock Oxford University Press, Oxford, 2002.

\bibitem{CampHern}
Felix Campelo and Aurora Hernandez-Machado.
\newblock Dynamic model and stationary shapes of fluid vesicles.
\newblock {\em Eur. Phys. Journal E-Soft Matter}, 20(1):37--45, 2006.

\bibitem{CampHern2}
Felix Campelo and Aurora Hernandez-Machado.
\newblock Shape instabilities in vesicles: A phase-field model.
\newblock {\em The European Physical Journal}, 143(1):101--108, 2007.

\bibitem{Can}
Peter~B. Canham.
\newblock The minimum energy of bending as a possible explanation of the
  biconcave shape of the human red blood cell.
\newblock {\em J. Theor. Biol.}, 26:61--81, 1970.

\bibitem{Chen}
Xinfu Chen.
\newblock Global asymptotic limit of solutions of the cahnÐhilliard equation.
\newblock {\em J. Differential Geometry}, 44(2):262--311, 1996.

\bibitem{DG}
Ennio De~Giorgi.
\newblock Some remarks on {$\Gamma$}-convergence and least squares method.
\newblock In {\em Composite media and homogenization theory (Trieste, 1990)},
  volume~5 of {\em Progr. Nonlinear Differential Equations Appl.}, pages
  135--142. Birkh\"auser Boston, Boston, MA, 1991.

\bibitem{DelOriVar}
Silvano Delladio.
\newblock Do generalized {G}auss graphs induce curvature varifolds?
\newblock {\em Boll. Un. Mat. Ital. B}, 10(4):991--1017, 1996.

\bibitem{Del}
Silvano Delladio.
\newblock Special generalized {G}auss graphs and their application to
  minimization of functionals involving curvatures.
\newblock {\em J. Reine Angew. Math.}, 486:17--43, 2007.

\bibitem{DuStokes}
Qiang Du, Chun Liu, and Manlin Li.
\newblock Analysis of a phase-field {N}avier-{S}tokes vesicle-fluid interaction
  model.
\newblock {\em Discrete Contin. Dyn. Syst.}, 8(3):539--556, 2007.

\bibitem{DuSpontCurv}
Qiang Du, Chun Liu, Rolf Ryham, and Xiaoqang Wang.
\newblock Modeling the spontaneous curvature effects in static cell membrane
  deformations by a phase field formulation.
\newblock {\em Commun. Pure Appl. Anal.}, 4(3):537--548, 2005.

\bibitem{DuWill}
Qiang Du, Chun Liu, Rolf Ryham, and Xiaoqiang Wang.
\newblock A phase field formulation of the willmore problem.
\newblock {\em Nonlinearity}, 18(3):1249--1267, 2005.

\bibitem{DuCa}
Qiang Du, Chun Liu, Rolf Ryham, and Xiaoqiang Wang.
\newblock Diffuse interface energies capturing the {E}uler number: relaxation
  and renormalization.
\newblock {\em Commun. Math. Sci.}, 8(1):233--242, 2007.

\bibitem{DuUno}
Qiang Du, Chun Liu, and Xiaoqiang Wang.
\newblock A phase field approach in the numerical study of the elastic bending
  energy for vesicle membranes.
\newblock {\em J. Comput. Phys.}, 198(2):450--468, 2004.

\bibitem{Dutopo}
Qiang Du, Chun Liu, and Xiaoqiang Wang.
\newblock Retrieving topological information for phase field models.
\newblock {\em SIAM J. Appl. Math.}, 65(6):1913--1932, 2005.

\bibitem{DuDue}
Qiang Du, Chun Liu, and Xiaoqiang Wang.
\newblock Simulating the deformation of vesicle membranes under elastic bending
  energy in three dimensions.
\newblock {\em J. Comput. Phys.}, 212(2):757--777, 2006.

\bibitem{Dumulti}
Qiang Du and Xiaoqiang Wang.
\newblock Modelling and simulations of multi-component lipid membranes and open
  membranes via diffuse interface approaches.
\newblock {\em J. Math. Biol.}, 56(3):347--371, 2008.

\bibitem{Eva}
Evan~A. Evans.
\newblock Bending resistance and chemically induced moments in membrane
  bilayers.
\newblock {\em Biophys. J.}, 14:921--931, 1974.

\bibitem{FarGar}
Hassan~M. Farshbaf-Shaker and Harald Garcke.
\newblock Thermodynamically consistent higher order phase field
  {N}avier-{S}tokes models with applications to biological membranes.
\newblock {\em Preprint}, 2009.

\bibitem{G}
Enrico Giusti.
\newblock {\em Minimal {S}urfaces and {F}unctions of {B}ounded {V}ariation}.
\newblock Birkh\"auser, Boston, 1984.

\bibitem{Gross}
Karsten Gro{\ss}e-Brauckmann.
\newblock New surfaces of constant mean curvature.
\newblock {\em Math. Z.}, 214:527--565, 1993.

\bibitem{GrossPhd}
Karsten Gro{\ss}e-Brauckmann.
\newblock {\em Complete embedded constant mean curvature surfaces}.
\newblock Habilitationsschrift. Universit\"at Bonn, Bonn, 1998.

\bibitem{Hu}
John~E. Hutchinson.
\newblock Second fundamental form for varifolds and the existence of surfaces
  minimising curvature.
\newblock {\em Indiana Univ. Math. J.}, 35(1):45--71, 1986.

\bibitem{Mant}
Carlo Mantegazza.
\newblock Curvature varifolds with boundary.
\newblock {\em J. Differential Geom.}, 43(4):807--843, 1996.

\bibitem{MM}
Luciano Modica and Stefano Mortola.
\newblock Un esempio di {$\Gamma$}-convergenza.
\newblock {\em Boll. Un. Mat. Ital. B (5)}, 14(1):285--299, 1977.

\bibitem{Mose01}
Roger Moser.
\newblock A generalization of {R}ellich's theorem and regularity of varifolds
  minimizing curvature, 2001.

\bibitem{MosDG}
Roger Moser.
\newblock A higher order asymptotic problem related to phase transitions.
\newblock {\em SIAM J. Math. Anal.}, 37(712--736):1--20, 2005.

\bibitem{PaTonz}
Pablo Padilla and Yoshihiro Tonegawa.
\newblock On the convergence of stable phase transitions.
\newblock {\em Comm. Pure App. Math.}, 51(6):551--579, 1998.

\bibitem{PellRoeg}
Mark~A. Pelletier and Matthias R\"oger.
\newblock Partial localization, lipid bilayers, and the elastica functional.
\newblock {\em Arch. Rational Mech. Anal.}, 193(3):475--537, 2009.

\bibitem{Petrov}
Alexander~G. Petrov.
\newblock {\em The Lyotropic State of Matter: Molecular Physics and Living
  Matter Physics}.
\newblock Gordon and Breach, Amsterdam, 1999.

\bibitem{RS}
Matthias R\"oger and Reiner Sch\"atzle.
\newblock On a modified conjecture of {D}e {G}iorgi.
\newblock {\em Math. Z.}, 254(4):675--714, 2006.

\bibitem{SchaJDG}
Reiner Sch\"atzle.
\newblock Hypersurfaces with mean curvature given an ambient sobolev function.
\newblock {\em J. Differential Geom.}, 58(3):371--420, 2001.

\bibitem{SiegKoz-Gauss}
David~P. Siegel and Michael~M. Kozlovy.
\newblock The {G}aussian {C}urvature {E}lastic {M}odulus of
  {N}-{M}onomethylated {D}ioleoylphosphatidylethanolamine: {R}elevance to
  {M}embrane fusion and lipid phase behavior.
\newblock {\em Biophysical Journal}, 87:366--374, 2004.

\bibitem{Si}
Leon Simon.
\newblock {\em Lectures on {G}eometric {M}easure {T}heory}, volume~3 of {\em
  Proceedings of the Centre for Mathematical Analysis, Australian National
  University}.
\newblock Australian National University Centre for Mathematical Analysis,
  Canberra, 1983.

\bibitem{TKS-Gauss}
Richard~H. Templer, Bee~J. Khoo, and John~M. Seddon.
\newblock Gaussian curvature modulus of an amphiphilic monolayer.
\newblock {\em Langmuir}, 14(26):7427--7434, 1998.

\bibitem{Tonesff}
Yoshihiro Tonegawa.
\newblock On stable critical points for a singular perturbation problem.
\newblock {\em Comm. Analysis and Geometry}, 13(2):439--459, 2005.

\bibitem{ToneNotes}
Yoshihiro Tonegawa.
\newblock Applications of geometric measure theory to two-phase separation
  problems.
\newblock {\em Sugaku Expositions}, 21(1):97--115, 2008.

\bibitem{ToneYuko}
Yoshihiro Tonegawa and Yuko Nagase.
\newblock A singular perturbation problem with integral curvature bound.
\newblock {\em Hiroshima Math. Journal}, 37(3):455--489, 2007.

\bibitem{Wango}
Xiaoqiang Wang.
\newblock Asymptotic analysis of phase field formulations of bending elasticity
  models.
\newblock {\em SIAM J. Math. Anal.}, 39(5):1367--1401, 2008.

\end{thebibliography}
\bibliographystyle{plain}
\end{document}